\documentclass[a4paper,12pt]{amsart}
\usepackage{amsthm,amsmath,amssymb}
\usepackage{pdfsync}
\usepackage{tikz-cd}
\usepackage{mathscinet}
\usepackage[notext]{stix2}
\def\join{\,\ast\,}
\def\joinB{\,\ast_{B}\,}
\def\Sq{\operatorname{Sq}}
\def\Hom{\operatorname{Hom}}
\date{9th February 2022} % from scratch
\date{\today}
\allowdisplaybreaks[4]
\theoremstyle{plain}
  \newtheorem{thm}{Theorem}[section]
  
  \newtheorem{prop}[thm]{Proposition}
  \newtheorem{exer}[thm]{Exercise}
  
  \newtheorem{prob}[thm]{Problem}
  \newtheorem{conj}[thm]{Conjecture}
\theoremstyle{definition}
  \newtheorem{defn}[thm]{Defninition}
  \newtheorem{expl}[thm]{Example}
\theoremstyle{remark}
  \newtheorem{rem}[thm]{Remark}
\def\hooklongrightarrow{\ensuremath{\lhook\joinrel\longrightarrow}}
\def\twoheadlongrightarrow{\ensuremath{\relbar\joinrel\twoheadrightarrow}}

\makeatletter
\@namedef{subjclassname@2010}{%
  \textup{2010} Mathematics Subject Classification}
\makeatother

\def\Img{\operatorname{Im}}
\def\Min{\operatorname{Min}}
\def\Max{\operatorname{Max}}
\def\cd{\operatorname{cd}}
\def\gd{\operatorname{gd}}
\def\crit#1{\operatorname{Crit}(#1)}
\def\ball#1{\operatorname{Ball}(#1)}
\def\gcat#1{\operatorname{{\it g}Cat}(#1)}
\def\wcat#1{\operatorname{{\it w}cat}(#1)}
\def\Imm#1{\operatorname{Imm}(#1)}
\def\double#1{#1 \!\times\! #1}
\def\paired#1#2{#1 \!\times\! #2}
\def\mathbold#1{\text{$\mathbb{#1}$}}
\makeatletter
\def\emptyarg{}
\def\@wgt#1{\def\thisarg{#1}\ifx\thisarg\emptyarg\operatorname{wgt}\else\operatorname{wgt}\hspace{.06em}(#1)\fi}
\def\@@wgt[#1]#2{\def\thisarg{#2}\ifx\thisarg\emptyarg\operatorname{wgt}_{#1}\else\operatorname{wgt}_{#1}(#2)\fi}
\def\wgt{\@ifnextchar[{\@@wgt}{\@wgt}}
\def\@Mwgt#1{\def\thisarg{#1}\ifx\thisarg\emptyarg\operatorname{Mwgt}\else\operatorname{Mwgt}\hspace{.06em}(#1)\fi}
\def\@@Mwgt[#1]#2{\def\thisarg{#2}\ifx\thisarg\emptyarg\operatorname{Mwgt}_{#1}\else\operatorname{Mwgt}_{#1}(#2)\fi}
\def\Mwgt{\@ifnextchar[{\@@Mwgt}{\@Mwgt}}
\def\@cuplen#1{\def\thisarg{#1}\ifx\thisarg\emptyarg\operatorname{cup}\else\operatorname{cup}\hspace{.06em}(#1)\fi}
\def\@@cuplen[#1]#2{\def\thisarg{#2}\ifx\thisarg\emptyarg\operatorname{cup}_{#1}\else\operatorname{cup}_{#1}(#2)\fi}
\def\cuplen{\@ifnextchar[{\@@cuplen}{\@cuplen}}
\def\@catlen#1{\def\thisarg{#1}\ifx\thisarg\emptyarg\operatorname{catlen}\else\operatorname{catlen}\hspace{.06em}(#1)\fi}
\def\@@catlen[#1]#2{\def\thisarg{#2}\ifx\thisarg\emptyarg\operatorname{catlen}_{#1}\else\operatorname{catlen}_{#1}(#2)\fi}
\def\catlen{\@ifnextchar[{\@@catlen}{\@catlen}}
\def\@Cat#1{\def\thisarg{#1}\ifx\thisarg\emptyarg\operatorname{Cat}\else\operatorname{Cat}\hspace{.06em}(#1)\fi}
\def\@@Cat[#1]#2{\def\thisarg{#2}\ifx\thisarg\emptyarg\operatorname{Cat}_{#1}\else\operatorname{Cat}_{#1}(#2)\fi}
\def\Cat{\@ifnextchar[{\@@Cat}{\@Cat}}
\def\@cat#1{\def\thisarg{#1}\ifx\thisarg\emptyarg\operatorname{cat}\else\operatorname{cat}\hspace{.06em}(#1)\fi}
\def\@@cat[#1]#2{\def\thisarg{#2}\ifx\thisarg\emptyarg\operatorname{cat}_{#1}\else\operatorname{cat}_{#1}(#2)\fi}
\def\cat{\@ifnextchar[{\@@cat}{\@cat}}
\def\@catFH#1{\def\thisarg{#1}\ifx\thisarg\emptyarg\operatorname{cat^{FH}}\else\operatorname{cat^{FH}}\hspace{.02em}(#1)\fi}
\def\@@catFH[#1]#2{\def\thisarg{#2}\ifx\thisarg\emptyarg\operatorname{cat^{FH}}_{#1}\else\operatorname{cat^{FH}}_{#1}(#2)\fi}
\def\catFH{\@ifnextchar[{\@@catFH}{\@catFH}}
\def\@catBG#1{\def\thisarg{#1}\ifx\thisarg\emptyarg\operatorname{cat^{BG}}\else\operatorname{cat^{BG}}\hspace{.02em}(#1)\fi}
\def\@@catBG[#1]#2{\def\thisarg{#2}\ifx\thisarg\emptyarg\operatorname{cat^{BG}}_{#1}\else\operatorname{cat^{BG}}_{#1}(#2)\fi}
\def\catBG{\@ifnextchar[{\@@catBG}{\@catBG}}
\def\@catB#1{\def\thisarg{#1}\ifx\thisarg\emptyarg\operatorname{cat}_{B}\else\operatorname{cat}_{B}(#1)\fi}
\def\@@catB#1{\def\thisarg{#1}\ifx\thisarg\emptyarg\operatorname{cat}_{B}\else\operatorname{cat}_{B}(#1)\fi}
\def\catB{\@ifnextchar[{\@@catB}{\@catB}}
\def\@catBB#1{\def\thisarg{#1}\ifx\thisarg\emptyarg\operatorname{cat}^{B}_{B}\else\operatorname{cat}^{B}_{B}(#1)\fi}
\def\@@catBB#1{\def\thisarg{#1}\ifx\thisarg\emptyarg\operatorname{cat}^{B}_{B}\else\operatorname{cat}^{B}_{B}(#1)\fi}
\def\catBB{\@ifnextchar[{\@@catBB}{\@catBB}}
\def\@tc#1{\def\thisarg{#1}\ifx\thisarg\emptyarg\operatorname{tc}\else\operatorname{tc}\hspace{.06em}(#1)\fi}
\def\@@tc{\operatorname{tc}}
\def\tc{\@ifnextchar({\@@tc}{\@tc}}
\def\@secat#1{\def\thisarg{#1}\ifx\thisarg\emptyarg\operatorname{genus}\else\operatorname{genus}\hspace{.06em}(#1)\fi}
\def\@@secat{\operatorname{genus}}
\def\secat{\@ifnextchar({\@@secat}{\@secat}}
\def\@Secat#1{\def\thisarg{#1}\ifx\thisarg\emptyarg\operatorname{Genus}\else\operatorname{Genus}\hspace{.06em}(#1)\fi}
\def\@@Secat{\operatorname{Genus}}
\def\Secat{\@ifnextchar({\@@Secat}{\@Secat}}
\def\@tcm#1{\def\thisarg{#1}\ifx\thisarg\emptyarg\operatorname{tc^{M}}\else\operatorname{tc^{M}}(#1)\fi}
\def\@@tcm{\operatorname{tc^{M}}}
\def\tcm{\@ifnextchar({\@@tcm}{\@tcm}}
\def\@TC#1{\def\thisarg{#1}\ifx\thisarg\emptyarg\operatorname{TC}\else\operatorname{TC}\hspace{.06em}(#1)\fi}
\def\@@TC{\operatorname{TC}}
\def\TC{\@ifnextchar({\@@TC}{\@TC}}
\def\@Path#1{\def\thisarg{#1}\ifx\thisarg\emptyarg\operatorname{\mathcal{P}}\else\operatorname{\mathcal{P}}\hspace{.06em}(#1)\fi}
\def\@@Path{\operatorname{\mathcal{P}}}
\def\Path{\@ifnextchar({\@@Path}{\@Path}}
\makeatother
\def\Conelen#1{\operatorname{Cl}(#1)}
\def\cwgt#1{\operatorname{cwgt}(#1)}
\def\fLoop#1{\operatorname{\mathcal{L}}\hskip.1em(#1)}
\def\Loop#1{\operatorname{\Omega}\hskip.1em(#1)}
\def\Susp#1{\operatorname{\Sigma}\hskip.1em(#1)}

\def\LoopBB#1{\operatorname{\Omega}^{B}_{B}(#1)}
\def\product{\text{\large$\textstyle\Pi$}}
\def\fatvee{\text{\large$\text{T}$}}

\def\homeo{\approx}
\def\widebar{\overline}
\def\sq{{\mathcal S\kern-0.2em q}}

\def\id{\operatorname{id}}
\def\proj{\operatorname{pr}}

\def\comp{\smash{\lower-.1ex\hbox{\scriptsize$\circ$\,}}}
\def\fracinline#1/#2{\mbox{\raise0.5ex\hbox{\footnotesize$#1$}{\hskip-.1em$/$\hskip-.1em}\raise-0.5ex\hbox{\footnotesize$#2$}}}
\def\field{{\mathbb F}}

\def\integral{{\mathbb Z}}

\def\real{{\mathbb R}}
\def\complex{{\mathbb C}}
\def\quaternion{{\mathbb H}}

\def\RP{\operatorname{{\real}P}}
\def\CP{\operatorname{{\complex}P}}
\def\HP{\operatorname{{\quaternion}P}}

\begin{document}

%%%%%

\baselineskip18pt

%%%%%

\title[L-S theory]{Lusternik-Schnirelmann theory \\to Topological Complexity\\from $A_{\infty}$-view point}

\author[N. Iwase]{Norio Iwase}
\address{Faculty of Mathematics,
Kyushu University,
Motooka 744, Fukuoka 819-0395, Japan}
\email{iwase@math.kyushu-u.ac.jp}

\begin{abstract}
We are trying to look over the Lusternik-Schnirelmann theory (L-S theory, for short) and the Topological Complexity (TC, for short) as a natural extension of the L-S theory.
In particular, we focus on the impact of the ideas originated from E.~Fadell and S.~Husseini on both theories.
More precisely, we see how their ideas on a category weight and a relative category drive the L-S theory and the TC.
\end{abstract}



\thanks{The author was supported in part by Grant-in-Aid for Scientific Research (S) \#17H06128 and by Exploratory Research \#18K18713 from Japan Society for the Promotion of Science.}

\allowdisplaybreaks
\maketitle

%%%%%%%%%%%%%%%%%%%%%%%%%%%%%%%%%%%%%%%%%%%%%%%%%%%%%%%%%%%
In this article, we work in the category of CW-complexes and maps between them, but we often restrict ourselves into the full subcategory of pointed CW-complexes. 
The pointed and unpointed theories are very close in the usual homotopy theory, but we find that they are far apart, if we discuss the higher associativity of an H-space fibrewise, or parametrized.
Before starting the main part, we discuss about the higher associativity of an H-space, which is closely related to both (fibrewise) L-S theory and TC. 

%%%%%%%%%%%%%%%%%%%%%%%%%%%%%%%%%%%%%%%%%%%%%%%%%%%%%%%%%%%
\section{Higher Associativity}

In 1941, H.~Hopf introduced, in \cite{MR4784}, a notion of an H-space, as a space with a multiplication with a homotopy unit.
The idea attracted many authors such as A.~Borel, I.~James, J.~F.~Adams, H.~Toda, W.~Browder, J.~Hubbuck, M. Sugawara, A. Zabrodsky, R.~Kane, J.~Lin, M.~Mimura and many others including the author. But an H-space such as the unit sphere of Octanions $S^{7}$ fails to satisfy (higher) homotopy associativity.

In 1957, M.~Sugawara gave a criteria for a space to be an H-space or a homotopy associative H-space in \cite{MR86303,MR97066}.
Refining extensively the idea employed by M.~Sugawara, J.~D.~Stasheff introduced a homotopy-theoretical version of the Milnor filtration \cite{MR77122,MR77932} as $A_{\infty}$-forms and $A_{\infty}$-structures in \cite{MR0158400} and a few years later in \cite{MR0270372} in a more sophisticated form, while the definitions in \cite{MR0158400} and in \cite{MR0270372} are slightly different:
In \cite{MR0158400}, an $A_{\infty}$-form requires a strict unit with higher coherency conditions, while in \cite{MR0270372} it requires only a homotopy unit as an H-space.
Because the induction step of the proof in \cite{MR0270372} claiming the equality of the two definitions stops at the stage $3$, it puzzled us for a long time, and finally is resolved in \cite{AX12115741}, and so the two definitions are not the same but \textit{equivalent, up to homotopy}.

One goal of $A_{\infty}$ theory is described as Theorem \ref{thm:stasheff-main} which is due to J.~D.~Stasheff.

See Appendix for explanations on Milnor filtrations of the classifying space of a topological group, definitions of $A_{\infty}$-forms and $A_{\infty}$-structures for both spaces and maps, and Ganea's fibre-cofibre constructions for connected spaces.

\medskip

%%%%%%%%%%%%%%%%%%%%%%%%%%%%%%%%%%%%%%%%%%%%%%%%%%%%%%%%%%%
\section{L-S theory}\label{sect:LStheory}

Let $M$ be a differentiable closed manifold and $f$ a smooth function on $M$.
Then, one basic question arises: how many critical points does $f$ have?
Let us denote by $\crit{f}$ the number of critical points of a smooth function $f$ on $M$, and by $\crit{M}$ the minimum among all the number $\crit{f}$ where $f$ runs over all smooth functions on $M$.
\par
In 1934, L.~Lusternik and L.~Schnirelmann, who worked on variation theory, gave a lower bound for $\crit{M}$ as their `category' number $\cat{M} \ge 1$, in other words, $\cat{M} \le \crit{M}$ in \cite{LS1934}.
However, among people working on L-S theory, it became popular to define $\cat{X}$ one less than the original:
the Lusternik-Schnirelmann (L-S, for short) category of a connected pointed space $X$, denoted by $\cat{X}$, is the smallest integer $n \!\ge\! 0$ such that $X$ is covered by $n{+}1$ categorical open subsets, where a subset $U \subset X$ is said to be categorical if the inclusion map $U \hookrightarrow X$ is homotopic to the constant map at the base point.
If no such integer exists, we write $\cat{X}=\infty$.
\begin{thm}[L.~Lusternik and L.~Schnirelmann]
Let $M$ be a manifold.
Then we have
$$\cat{M} \le \crit{M}\!-\!1 \le \dim{M}.$$
\end{thm}

In 1939-41, R.~H.~Fox introduced, in \cite{MR2937306,MR4108}, a series of new ideas related to L-S theory. 
Firstly, for a manifold $M$, a ball category is denoted by $\ball{M}$ as the smallest integer $n \!\ge\! 0$ such that $M$ is covered by $n{+}1$ open balls in $M$.
Secondly, for a connected pointed space $X$, his version of strong category or a geometric category is denoted by $\gcat{X}$, as the smallest integer $n \!\ge\! 0$ such that $X$ is covered by $n{+}1$ contractible open subsets, which is later reformulated by T.~Ganea as his homotopy-invariant version of strong category.
We remark that $\ball{M}$ as well as $\gcat{X}$ are not homotopy invariant (see \cite{MR1990857}).
Thirdly, for a connected pointed space $X$, a finite sequence $\{A_{0},\ldots,A_{k}\}$ of closed subsets of $X$ a categorical sequence of length $k$, if $\{\ast\} \subset A_{0} \subset A_{1} \subset \cdots \subset A_{k-1} \subset A_{k}{=}X$ such that $\{A_{0},A_{1}\!\smallsetminus\!A_{0},\ldots,A_{k}\!\smallsetminus\!A_{k-1}\}$ are all categorical in $X$.
Then $\catlen{X}$ is the minimal length of all such sequences, which can also be reformulated (see \cite{MR2503529}) using Fadell-Husseini relative L-S category which will be introduced later.
\begin{thm}[R.~H.~Fox]
Let $M$ be a manifold, and $X$ a space.
Then we have
\begin{enumerate}
\item
$\cat{M} \le \gcat{M} \le \ball{M} \le \crit{M}{-}1 \le \dim{M}$, and
\item
$\cat{X}=\catlen{X} \le \gcat{X} \le \dim{X}$.
\end{enumerate}
\end{thm}

To determine L-S category, we need some computable lower bound.
Let $h$ be a multiplicative cohomology theory.
Classically, the cup-length of $X$ w.r.t. $h$, denoted by $\cuplen{X;h}$, is the supremum of $k \!\ge\! 0$ such that there exists a non-zero $k$-fold product in the reduced theory $\tilde{h}^{\ast}(X)$.
When $h$ is the ordinary cohomology with coefficients in $R$ a ring with unit, then we denote $\cuplen{X;R} = \cuplen{X;h}$.
The following theorem is well known, while $\cuplen{X;h}$ is much less than $\dim{X}$ in a number of examples.

\begin{thm}
\label{thm:inequality}
$\cuplen{X;h} \le \cat{X} \le \gcat{X} \le \dim{X}$.
\end{thm}

In 1960, I.~Berstein and P.~J.~Hilton \cite{MR126276} gave a criterion for $\cat{C_{f}}=2$ where $C_{f}$ is the mapping cone of a map $f : X \to Y$, in terms of their version of a Hopf invariant $H_1(f) \in [\Sigma X, \Omega \Sigma Y{\ast} \Omega \Sigma Y]$ for a map $f : \Sigma X \to \Sigma Y$, where $A{*}B$ denotes the join of spaces $A$ and $B$.
In addition, its higher version $H_{m}(f)$ is used to disprove the Ganea conjecture on L-S category (see \cite{MR1642747,MR1905835}). 
They also introduced a new lower bound called weak category, denoted by $\wcat{X}$, which is the supremum of $k \!\ge\! 0$ such that the reduced iterated diagonal $\overline{\Delta}^{k} : X \to \wedge^{k}X$ is non-trivial.
Similarly, we denote by $\cuplen{X}$ the supremum of $k \!\ge\! 0$ such that the reduced iterated diagonal $\widebar{\Delta}^{k} : X \to \wedge^{k}X$ is stably non-trivial.
For example, M.~Mimura and the author showed in \cite{MR2039767} that $\widebar{\Delta}^{n+2} : \mathrm{Sp}(n) \to \wedge^{n+2}\mathrm{Sp}(n)$ is non-trivial and $\cat{\mathrm{Sp}(n)} \!\ge\! n{+}2$ for $n \!\ge\! 3$ by showing that $\widebar{\Delta}^{4}$ $:$ $\mathrm{Sp}(3) \overset{q}\twoheadrightarrow S^{18} \xrightarrow{\nu^{2}} S^{12} \overset{j}\hookrightarrow \wedge^{4}\mathrm{Sp}(3)$ is non-trivial, using Toda secondary composition \cite{MR0143217} which is obtained independently by L.~Fern\'andez-Su\'arez, A.~G\'omez-Tato, J.~Strom and D.~Tanr\'e in \cite{MR2022385}. 

\begin{thm}
\label{thm:weakcat}
$\cuplen{X;h} \le \cuplen{X} \le \wcat{X} \le \cat{X}$.
\end{thm}

In 1962, I.~Berstein and T.~Ganea defined in \cite{MR139168} a L-S category of a map.
For a map $f$ from a space $K$ to a connected pointed space $X$, they defined $\cat{f}$ as the smallest integer $n \!\ge\! 0$ such that $K$ is covered by $n{+}1$ open subsets on which the restriction of $f$ is homotopic to the constant map at the base point.
If no such integer exists, we write $\cat{f}=\infty$.
For an inclusion $i : K \hookrightarrow X$, let us denote $\cat{X;K}=\cat{i}$.

In 1963, J.~D.~Stasheff introduced a notion of an $A_{m}$-space using $A_{m}$-forms and $A_{m}$-structures to amount the higher associativity for an H-space in \cite{MR0158400}.
The existence of an $A_{m}$-structure implies the `standard' $A_{m}$-structure which enjoys the following three properties (see \cite{MR0158400}).
The advantage to consider Stasheff's construction is that the cell structure of his `standard' $A_{m}$-structure is quite understandable.
\begin{align}&
P^{0}(G) =\{\ast\},\tag{Z}
\\&
G \hookrightarrow E^{n}(G) \twoheadrightarrow P^{n-1}(G) \ \ \text{is a fibre sequence, $m \!\ge\! n \!\ge\! 1$,}\tag{F'}
\\&
E^{n}(G) \twoheadrightarrow P^{n-1}(G) \hookrightarrow P^{n}(G)\ \ \text{is a cofibre sequence, $m \!\ge\! n \!\ge\! 1$.}\tag{C}
\end{align}
When $m\!=\!\infty$, they enjoys the following three properties (see \cite{MR0158400,MR0270372}, for details).
\begin{align}&
P^{0}(G) =\{\ast\},\tag{Z}
\\&
E^{n}(G) \twoheadrightarrow P^{n-1}(G) \hookrightarrow P^{\infty}(G)
\ \ \text{is a fibre sequence up to homotopy, $n \!\ge\! 1$,}\tag{F}
\\&
E^{n}(G) \twoheadrightarrow P^{n-1}(G) \hookrightarrow P^{n}(G)\ \ \text{is a cofibre sequence, $n \!\ge\! 1$.}\tag{C}
\end{align}

Using $A_{\infty}$-maps between $A_{\infty}$-spaces (see \S \ref{app-sect:higherassociativity}), J.~D.~Stasheff showed in \cite{MR0158400,MR0270372} the following fundamental result saying that his $P^{\infty}$ gives the inverse functor to $\Omega$.
\begin{thm}[J.~D.~Stasheff]\label{thm:stasheff-main}
For a connected CW complex $X$, $\Loop{X}$ is an $A_{\infty}$-space and $P^{\infty}\Loop{X} \simeq X$.
For an $A_{\infty}$-space $G$, which has the homotopy type of a connected CW complex, $\Loop{P^{\infty}G}$ is $A_{\infty}$-homotopy equivalent to $G$.
\end{thm}
From now on, $\{(E^{n+1},P^{n},p_{n+1},i_{n},e_{n})\}_{n \ge 0}$ consisting of spaces $E^{n+1}$ and $P^{n}$, and maps $p_{n+1} : E^{n+1} \to P^{n}$, $i_{n} : P^{n} \to P^{n+1}$ and $e_{n} : P^{n} \to X$ for all $n \!\ge\! 0$ is called a fibre-cofibre construction for a pointed space $X$, if it satisfies the following three conditions.
\begin{align}&
P^{0} =\{\ast\},\tag{Z}
\\&
E^{n} \xrightarrow{p_{n}} P^{n-1} \xrightarrow{e_{n\!-\!1}} X\ \ \text{is a fibre sequence up to homotopy for $n \ge 1$,}\tag{F}
\\&
E^{n} \xrightarrow{p_{n}} P^{n-1} \xrightarrow{i_{n\!-\!1}} P^{n}\ \ \text{is a cofibre sequence up to homotopy for $n \ge 1$.}\tag{C}
\end{align}

In 1967, T.~Ganea introduced in \cite{MR229240} a strong category $\Cat{X}$ as a homotopy invariant version of a Fox strong category or a geometric category $\gcat{X}$. 
More precisely, $\Cat{X}$ is the minimum among all the number $\gcat{Y} \!\ge\! 0$ where $Y$ runs over all spaces with the same homotopy type of $X$.
T.~Ganea also showed that $\Cat{X}$ is characterized as follows:
for a connected pointed space $X$, $\Cat{X}$ is $0$ if $X$ is contractible and, otherwise, is equal to the smallest integer $m$ such that there is a series of cofibre sequences
$\{K_{i} \to F_{i-1}\hookrightarrow F_{i} \;|\; 1 \!\le\! i \!\le\! m\}$ with $F_{0}= \{\ast\}$ and $F_{m}\simeq X$.
Such cofibre sequence is called a cone-decomposition of $X$, and the smallest integer $m$ is called the cone-length of $X$ denoted here by $\Conelen{X}$, which is a similar idea to Fox's categorical length (see Fox \cite{MR2937306,MR4108}) with category replaced by Ganea's strong category.

\begin{thm}[T.~Ganea]
\label{thm:inequality2}
$\Cat{X}\!-\!1 \!\le\! \cat{X} \!\le\! \Cat{X}\!=\!\Conelen{X} \!\le\! \gcat{X}$.
\end{thm}
We know that $\cat{X} \le 1$ if and only if $X$ is a co-H-space, while $\Cat{X} \le 1$ if and only if $X$ is desuspendable, that is, $X=\Susp{Y}$ for some $Y$.
Since there are lots of co-H-spaces which are not desuspendable, $\cat{}$ and $\Cat{}$ are actually different, in general.
For instance, a series of counter examples to the Ganea conjecture for a co-H-spaces in \cite{MR1808218} gives such examples.
This reminds us a result on discrete groups by S.~Eilenberg and T.~Ganea \cite{MR85510}.

\begin{thm}[S.~Eilenberg and T.~Ganea]
\label{thm:inequality3}
Let $G$ be a discrete group. Then we have
$\cd{G} \!\le\! \gd{G}$, where the equality holds except for the case when $(\gd{G},\cd{G})=(3,2)$.
\end{thm}
Because there is a space $X$ with $\Cat{X}=\cat{X}\!+\!1$, there could be a group $G$ with $\gd{G}=\cd{G}\!+\!1=3$.
But at this moment, the author doesn't have any idea on it.

\medskip

In 1967, T.~Ganea gave a fibre-cofibre construction $\{(E_{n},G_{n},p_{n},i_{n-1},e_{n})\}$ for a space $X$ as \cite[Proposition 2.2]{MR229240}, which is often called `Ganea construction'.

\begin{prop}[T.~Ganea]\label{ganea-fibrecofibre}
$\cat{X} \le m$ if and only if $p_{n}$ has a right homotopy inverse.
\end{prop}
To examine the existence of the right homotopy inverse of $p_{n}$ in the above proposition seriously, it might be useful to give the cell structure of the Ganea construction.

\medskip

In 1966-1976, an interesting result was obtained independently by A.~S.~\v{S}varc \cite{AMST1966Svarc} and by I.~Berstein \cite{MR400212} as Theorem \ref{thm:berstein-svarc} below in the case when $n \!\ge\! 3$.
The result was based on the study of I.~Berstein \cite{MR96221} in 1958 introducing the idea of cup-length of the cohomology theory with local coefficients. 
In 2009, A.~Dranishnikov and Y.~Rudyak \cite{MR2475974} gave a new proof of the result to include the case when $n\!=\!2$.

For a cell complex $X$, we denote the universal covering space of $X$ by $p : \widetilde{X} \to X$ and $\widetilde{X}_{0}=p^{-1}(\ast) \homeo \pi=\pi_{1}(X)$.
Then the group ring $\integral\pi$ acts on $C_{*}(\widetilde{X},\widetilde{X}_{0})$, and so we have a cochain group $C^{*}(X,\ast;M)=\Hom_{\integral\pi}(C_{*}(\widetilde{X},\widetilde{X}_{0}),M)$ of $\integral\pi$-homomorphisms for $\integral\pi$-module $M$.
Let $I(\pi)$ be the kernel of the augmentation $\varepsilon : \integral\pi \to \integral$.
Then we have the `fundamental class' $\bar{\frak{b}}$ $\in$ $H^{1}(X,\ast;I(\pi)) \cong \Hom(H_{1}(\widetilde{X},\widetilde{X}_{0}),I(\pi)) \cong \Hom(I(\pi),I(\pi))$ corresponding to the identity homomorphism.
The element $\frak{b}=j^{*}\bar{\frak{b}} \in H^{1}(X,I(\pi))$ is often called the Berstein-\v{S}varc class, where $j : X \hookrightarrow (X,\ast)$ is the canonical inclusion.
\begin{thm}[I.~Berstein, A.~S.~\v{S}varc, A.~Dranishnikov-Y.~Rudyak]\label{thm:berstein-svarc}
For a space satisfying $\cat{X}=\dim{X}=n \ge 2$, $\frak{b}^{n} \not= 0$ in $H^{n}(X,I(\pi)^{n})$, where $I(\pi)^{n}=I(\pi) \otimes \cdots \otimes I(\pi)$.
\end{thm}

While the definition of L-S category is fairly simple, it is not quite easy to determine the number of a given space.
As for an upper bound, strong category or cone-length gives a candidate for $\cat{X}$.
As for a lower bound, the cup-length gives a candidate for $\cat{X}$.
If the two candidates are fortunately the same, we can conclude that is the one we want.
But in general, the two candidates are far apart from each other.

In 1975-76, W.~Singhof \cite{MR391075,MR425952} determined the L-S category of a series of Lie groups by using cup-length arguments together with an open ball covering as $\cat{\mathrm{U}(n)}=n$ and $\cat{\mathrm{SU}(n)}=n{-}1$.
But unfortunately, this method can not be applied to other types of compact connected Lie groups, because we do not have general machinery for obtaining an open ball covering whose cardinality is one more than the cup-length.

In 1992, E.~Fadell and S.~Husseini \cite{MR1317569} introduced a new lower bound for $\cat{X}$ called `category weight' defined for an element in the reduced cohomology theory of $X$ depending on the topological structure of $X$ and the subspace category in $X$.
That means the invariant is topological but not homotopical.
Let $h$ be a cohomology theory.
\begin{defn}[E.~Fadell and S.~Husseini]
\label{def:originalcategoryweight}
$$
\cwgt{u} = \Min\{ m \ge 0 \mid \forall\,A \underset{\text{closed}}\subset X \ u|_{A}=0 \ \text{if} \ \cat{X;A}<m \},\quad u \in \tilde{h}^{*}(X),
$$
where $\cat{X;A}=\catBG{X,A}$ denotes the L-S category of $A$ in $X$, which is the same as $\cat{i}$ the L-S category of the inclusion map $i : A \hookrightarrow X$.
\end{defn}

The homotopy invariant version of `category weight' is introduced independently by J.~Strom and Y.~Rudyak, and is applied to determine L-S category of various spaces (see Y.~Rudyak and J.~Oprea \cite{MR1686579}, for example) and to solve Arnold conjecture (see Y.~Rudyak \cite{MR1686583} and/or K.~Fukaya and K.~Ono \cite{MR1688434}).

\begin{defn}[Y.~Rudyak \cite{MR1679849} and J.~Strom \cite{MR2696593}]
\label{defn:categoryweight}~
\begin{enumerate}
\item
$\wgt{u;h} = \Min\{m \ge 0 \mid f^{*}(u)\!=\!0 \ \text{if} \ \cat{f} \!<\! m\}$, \ \ $u \in \tilde{h}^{*}(X)$,
\item
$\wgt{X;h} = \Max\{\wgt{u;h} \mid 0 \!\not=\! u \!\in\! \tilde{h}^{*}(X)\}$.
\end{enumerate}
\end{defn}

Let us denote by $\wgt{X}$ the maximum of $\wgt{X;h}$ where $h$ runs over all cohomology theories.
Then we have the following.

\begin{thm}
\label{thm:inequality4}
$\wgt{X;h} \le \wgt{X} \le \cat{X}$. 
\end{thm}

In 1994, E.~Fadell and S.~Husseini \cite{MR1397466} introduced a relative version of L-S category for a pair. 
The relative L-S category of a pair $(K,L)$, denoted by $\cat{K,L}=\catFH{K,L}$, is the smallest integer $n \!\ge\! 0$ such that $K$ is covered by $n{+}1$ open subsets $U_{1}, \ldots, U_{n}$ and $V$, where $U_{i}$\,s are categorical in $K$ and $V$ is compressible relative to $L$ into $L$.
If no such integer exists, we write $\cat{K,L}=\infty$.
When $X \!\supset\! K \!\supset\! L \!\supset\! A$, it might be useful to define $\cat{X\,;K,L\!:\!A}$ as the smallest integer $n \!\ge\! 0$ such that $K$ is covered by $n{+}1$ open subsets $U_{1}, \ldots, U_{n}$ and $V$, where $U_{i}$\,s are compressible in $X$ relative to $A$ into $A$ and $V$ is compressible in $X$ relative to $A$ into $L$ (compare with \cite{MR2503529}).
This idea characterizes the Fox categorical sequence to obtain an upper bound for $\cat{X}$ using a cone-decomposition enhanced with higher Hopf invariants.

In 1994, O.~Cornea showed in \cite{MR1259517} that $\cat{P^{m}\Loop{X}} = \Min\{m, \cat{X}\}$.

In 1998, L-S category is described in terms of an $A_{\infty}$-structure of the loop space of a well-pointed space or a space with non-degenerate base point, using the Whitehead definition of L-S category, which is performed in \cite{MR1642747}.
The Whitehead definition says $\cat{X} \le m$ if and only if the $m{+}1$-fold diagonal $\Delta^{m+1} : X \to \overset{m+1}\product X$ is compressible into the fat wedge $\overset{m+1}\fatvee X$, where $\overset{k+1}\fatvee X$, $1 \!\le\! k \!\le\! m$ is defined by induction on $k$ as follows:
\begin{align*}&
(\overset{1}{\product}X,\overset{1}{\fatvee}X) = (X,\ast)\quad\text{and}\quad 
(\overset{k+1}{\product}X,\overset{k+1}{\fatvee}X) = (\overset{k}{\product}X \times X,\overset{k}{\fatvee}X \times X \cup \overset{k}{\product}X \times\ast).
\end{align*}
The homotopy pull-back of $\overset{m+1}{\fatvee}X \hooklongrightarrow \overset{m+1}{\product}X$ and $X \xrightarrow{\Delta^{m+1}} \overset{m+1}{\product}X$ has the homotopy type of the projective $m$-space $P^{m}\Loop{X}$ of $\Loop{X}$ and the homotopy pull-back of $\overset{m+1}{\fatvee}X \hooklongrightarrow \overset{m+1}{\product}X$ and $\ast \longrightarrow \overset{m+1}{\product}X$ has the homotopy type of $E^{m+1}\Loop{X}$ the $m{+}1$-fold join of $\Loop{X}$, by using the cube lemma (see \cite{MR402694}).

\medskip

Here, let us emphasize the difference between two notions -- a homotopy equivalence and an $A_{\infty}$-equivalence, because these notions are crucial ideas to understand the $A_{\infty}$-method on (fibrewise) L-S theory:
homotopy equivalent $A_{\infty}$-spaces do not always have the same $A_{\infty}$-structure.
For instance, $\mathrm{U}(n)$ and $\mathrm{SU}(n) \!\times\! S^{1}$ are homeomorphic Lie groups with homotopically different classifying spaces.
In other words, for a space $X=B\mathrm{SU}(n) \!\times {\complex}P^{\infty}$ and a Lie group $\mathrm{U}(n)$ of the same homotopy type as the monoid $\mathrm{SU}(n) \!\times\! S^{1} = \Omega{X}$, there does not exist a homotopy equivalence from $B\mathrm{U}(n)$ to $X$.
So, \cite[Exercise 2.16]{MR1990857} would be better to be replaced with something similar to the following.
\begin{exer}\label{exer:ganeaspaces}
Let $X$ and $Y$ be pointed spaces, and let
$\{(E^{n+1},P^{n},p_{n+1},i_{n},e_{n})\}_{n \ge 0}$ and $\{(F^{n+1},Q_{n},q_{n+1},j_{n},f_{n})\}_{n \ge 0}$ be their fibre-cofibre constructions, respectively.
Then a homotopy equivalence $\phi : X \to Y$ induces homotopy equivalences $\hat\Phi_{n} : E^{n} \to F^{n}$ and $\Phi_{n} : P^{n} \to Q^{n}$, $n \!\ge\! 0$ satisfying $\Phi_{n}{\comp}p_{n+1}=q_{n+1}{\comp}\hat\Phi_{n+1}$, $\Phi_{n+1}{\comp}i_{n}=j_{n}{\comp}\Phi_{n+1}$ and $f_{n}{\comp}\Phi_{n}=e_{n}$ for $n \!\ge\! 0$.
\end{exer}
As is claimed in \cite{MR1990857} for \cite[Exercise 2.16]{MR1990857}, it is not very hard to give a proof to the above exercise, and so we leave it to the reader.

Using Exercise \ref{exer:ganeaspaces}, the standard $A_{\infty}$-structure in \cite{MR0158400} for an $A_{\infty}$-space as well as the Milnor filtration of the classifying space for a topological group gives the fiber-cofibre construction.
Or, maybe we could say that a Ganea space in \cite{MR229240} of a space $X$ is nothing but a projective space of $\Loop{X}$ by the definition of an $A_{\infty}$-structure in \cite{MR0158400}.

\begin{thm}\label{ganea-iwase}
$\cat{X} \le m$ if and only if the natural map $e^{X}_{m} : P^{m}\Loop{X} \to X$ has a right homotopy inverse $\sigma$ in the category of pointed spaces.
\end{thm}

In 2007, a stronger version of category weight was introduced to determine the L-S category of $\mathrm{Spin}(9)$ (see  \cite{MR2272137,MR2331629}).
Let $\mathcal{A}_{p}$ be the modulo $p$ Steenrod algebra.
Then, since the map $e^{X}_{m}$ and $\sigma$ in Theorem \ref{ganea-iwase} gives an $\mathcal{A}_{p}$-homomorphism between $\mathcal{A}_{p}$-modules, the image of $(e^{X}_{m})^{*} : H^{*}(X;\field_{p}) \to H^{*}(P^{m}\Loop{X};\field_{p})$ is a direct summand of $H^{*}(P^{m}\Loop{X};\field_{p})$ as an $\mathcal{A}_{p}$-module.
Thus $\cat{X} \le m$ implies that the image of $(e^{X}_{m})^{*}$ is a direct summand of $H^{*}(P^{m}\Loop{X};\field_{p})$ as an $\mathcal{A}_{p}$-module.
Let $\Gamma$ be a subset of $h^{*}h$ the algebra of cohomology operations on $h$.
A direct summand of an $h^{*}h$-module is called $\Gamma$-direct summand, if it is closed under the action of $\Gamma \subset h^{*}h$.

\begin{defn}
\label{defn:categoryweight2}
For a space $X$ and a subset $\Gamma \subset h^{*}h$, we define
$$\Mwgt{X;\Gamma} = \Min\{ m \!\ge\! 0 \mid \text{$\Img(e^{X}_{m})^{*}$ is a $\Gamma$-direct summand}\}.$$
\end{defn}
For example, it is shown in \cite{MR2272137} that $\cat{\mathrm{Spin}(9)}$ $=$ $\Mwgt{\mathrm{Spin}(9);\{\Sq^{2}\}}$ $=$ $8$, while $\Mwgt{\mathrm{Spin}(9);\field_{2}}$ $=$ $\wgt{\mathrm{Spin}(9);\field_{2}}$ $=$ $6$.
Let us denote by $\Mwgt{X}$ the maximum of $\Mwgt{X;h^{*}h}$ where $h$ runs over all cohomology theories. Then we have the following theorem.

\begin{thm}
$\cuplen{X }\le \wgt{X} \le \Mwgt{X} \le \cat{X}$.
\end{thm}

It seems that we need a different idea to solve the following long-standing conjecture.

\begin{conj}
\label{conj:so}
$\cat{\mathrm{SO}(n)} = \cuplen{\mathrm{SO}(n)}$ for all $n \!\ge\! 1$.
\end{conj}
The conjecture is verified up to $n\!=\!10$, at this moment (see \cite{MR3520436}).

%%%%%%%%%%%%%%%%%%%%%%%%%%%%%%%%%%%%%%%%%%%%%%%%%%%%%%%%%%%
\section{Sectional Category and Topological Complexity}

We all know that the L-S category is representing a kind of complexity of a given space.
From the view point of Robot Motion Planning, we may think it is actually the same as the complexity for a Robot Motion Planning with a fixed Robot station.

Let us first recall the definition of a genus (see A.~S.~\v{S}varc \cite{MR0154284,MR0151982}) 
or a sectional category (see I.~M.~James \cite{MR516214}): In 1959-62, A.~S.~\v{S}varc introduced a homotopy invariant for a fibration $\varpi : P \to W$.
Following I.~M.~James \cite{MR516214}, we denote by $\Secat{\varpi}$ the smallest integer $n \!\ge\! 1$ such that $W$ is covered by $n$ sectional open subsets; that is, open subsets over which there is a section of $\varpi$.
For a map $f : Y \to W$, the sectional category of $f$ is defined as $\Secat{f}=\Secat{f^{*}\varpi}$ the sectional category of $f^{*}\varpi: f^{*}P \to Y$ the fibration induced from $\varpi : P \to W$ by $f$.

For two fibrewise spaces $p : E \to X$ and $q : F \to X$ over a space $X$, we can topologise $E \joinB F = \underset{x \in X}\bigcup\,(E_{x} \join F_{x})$ as a subspace of $E \join F$ (see I.~M.~James \cite{MR516214}) with a natural projection $E \joinB F \to X$.
Then by taking an iterated fibrewise join construction, we obtain $J^{m}(p) : J^{m}(E)=E \joinB \cdots \joinB E \to X$ ($m$-fold fibrewise join)

\begin{thm}[A.~S.~\v{S}varc]
$\Secat{\varpi} \le m$ if and only if the $m$-fold fibrewise join $J^{m}(\varpi) : J^{m}(P) \to W$ of $\varpi : P \to W$ admits a section.
\end{thm}

In 2003, M.~Farber \cite{MR1957228} introduced an idea of Topological Complexity for Robot Motion Planning (without a Robot Station).
Let $\Path{X}$ be the space of all paths on $X$ and $\pi : \Path{X} \to \double{X}$ the projection designating initial and terminal points, which is often called a Serre path fibration.
Then the original Topological Complexity (TC, for short) of a connected space $X$, denoted by $\TC{X}$, is defined to be $\Secat{\pi}$.

Just as in L-S category, people working on TC (including M.~Farber) have come to define TC as one less than the original definition as the minimum number of open sets with sections.
To avoid confusion, we denote $\secat{\varpi}=\Secat{\varpi}{-}1$ and $\tc{X}=\TC{X}-1$ for the new definition.
The following result by M.~Farber gave a strong impact not only to algebraic topologists but also to differential topologists.
\begin{thm}[M.~Farber \cite{MR1957228}]
For all $n \!\ge\! 1$, we have
$$\tc{\RP^{n}}=\Imm{\RP^{n}}-\delta_{n},\quad \delta_{n}=\left\{\vphantom{\bigg(}\!\!\!\!\right.\begin{array}{ll}1,&n=1,3,7\\0,&\text{otherwise,}\end{array}
$$
where $\Imm{M}$ is the Euclidean immersion dimension of a closed manifold $M$.
\end{thm}
This motivated many people to study TC more seriously.
Also the theoretical similarity with L-S theory drew attention of people working on L-S theory including the author.

In 2008, M.~Farber and M.~Grant \cite{MR2407101} introduced a TC version of `category weight' called TC-weight, and show its strength determining TC.
\begin{defn}[M.~Farber and M.~Grant]
Let $h$ be a cohomology theory and $Z^{*}_{\varpi} = \ker\{\Delta^{*} : h^{*}(\double{X}) \to h^{*}(X)\}$ the zero-divisors ideal.
\begin{enumerate}
\item
$\wgt[\varpi]{u;h} = \Min\{m \!\ge\! 0 \mid f^{*}(u)\!=\!0 \ \text{if} \ \secat{f} \!<\! m\}$, \ \ $u \!\in\! Z^{*}_{\varpi}$,
\item
$\wgt[\varpi]{X;h} = \Max\{\wgt[\varpi]{u;h} \mid 0 \!\not=\! u \!\in\! Z^{*}_{\varpi}\}$.
\end{enumerate}
\end{defn}

The original motivation of Robot Motion Planning leads us to consider various kinds of configuration spaces where robots are moving around. We can also consider a robot of a special shape which forces us to consider some advanced featured TCs including sequential or higher versions of TC.
More recently, leading researchers of TC are shifting to a more applied-side, introducing a parametrized TC (which takes account of obstacles that robots must avoid) and etc.

As for the original TC, there looks still a lot of things to learn from E.~Fadell and S.~Husseini.
In particular, we want to know about a relative TC and TC module weight.

%%%%%%%%%%%%%%%%%%%%%%%%%%%%%%%%%%%%%%%%%%%%%%%%%%%%%%%%%%%
\section{Fibrewise L-S theory}

In 2010, M.~Sakai and the author learned a fibrewise L-S theory from I.~M.~James and J.~R.~Morris \cite{MR1130605} and M.~C.~Crabb and I.~M.~James \cite{MR1646248}, and found that a fibrewise L-S theory and TC are almost the same.
The difference lies in the motion of a base point: James's fibrewise L-S category $\catBB{}$ corresponds to a `pointed' version of TC which we named `monoidal' TC denoted by $\tcm{}$, requiring that a robot does not move in the case when the initial and terminal states are the same.  The original TC corresponds to an `unpointed' version of a fibrewise L-S category which we named as fibrewise unpointed L-S category \cite{MR2556074,MR2923451}.
In other words, the original TC does allow a robot to move around, even if the given initial and terminal states of the robot are the same.

As is seen in the original paper by E.~Fadell and S.~Husseini \cite{MR1317569}, there is a strong connection between the category weight and the bar resolution of a group.
In fact, \cite{MR1905835,MR2331629}  showed the connection can be applied to every topological space.
Going into fibrewise, we see the same idea can be applied to every fibrewise pointed space.

A fibrewise space is a surjection $p : E \to X$ and is denoted by $E=(E,p,X)$, where $E$, $X$, $p$ are called the total space, the base space, and the projection, respectively.
Also a fibrewise pointed space is a pair of a fibrewise space $E=(E,p,X)$ and a cross-section $s$ of $p$, which is often denoted by $(E,s)=(E,p,X,s)$.
For a given fibrewise pointed space $(E,s)$, its fibrewise loop space $\LoopBB{E,s}$ can be constructed naturally in \cite{MR1130605}, which has a natural fibrewise $A_{\infty}$-form by \cite{MR2665235}.
A section $s : X \hookrightarrow E$ for $E=(E,p,X)$ is said to be a non-degenerate fibrewise base point, if it is a closed cofibration.
A fibrewise space with a non-degenerate fibrewise base point is called a fibrewise well-pointed space.

For example, the product space $E=\paired{X}{Y}$ is a fibrewise space $E=(E,\proj_{1},X)$ over $X$. Further, for any given map $f : X \to Y$, we obtain a cross-section $s_{f} : X \to E$ to $\proj_{1}$ by $s_{f}=(\paired{\id}{\,f})\comp\Delta : X \to \paired{X}{Y}$.
On the one hand, if $Y$ has a non-degenerate base point $\ast$, then $s_{*}$ is a non-degenerate fibrewise base point for $E$ where $* : X \to \{\ast\} \subset Y$ denotes the trivial map.
In this case, the associated fibrewise loop space of $E$ is nothing but the product space $X \!\times\! \Loop{Y}$.
On the other hand, if $(Y,X)$ is a CW-pair with closed cofibration $i : X \hookrightarrow Y$, then $(\paired{X}{Y},\Delta{X})$ is an NDR-pair and $s_{i}=(\paired{\id}{\,i})\comp\Delta : X \to \paired{X}{Y}$ is a non-degenerate fibrewise base point for $E$.
In this case, the associated fibrewise loop space of $E$ is nothing but the restriction to $X$ of the free loop space $\fLoop{Y}$ over $Y$.

The Whitehead definition of an L-S category for a well-pointed space can be extended to that for a fibrewise well-pointed space.
The fibrewise fat wedge $\overset{m+1}\fatvee_{\!\!\!\!B} E$ of a fibrewise pointed space $E=(E,s)$ over $X$ is given by induction on $k \!\ge\! 1$ by $(\overset{1}{\product}_{B}E,\overset{1}{\fatvee}_{\!B}E)$ $=$ $(E,s(X))$ and
\begin{align*}&
(\overset{k+1}{\product}_{\!\!B}E,\overset{k+1}{\fatvee}_{\!\!\!B}E) = (\overset{k}{\product}_{B}E \times_{B}E,\overset{k}{\fatvee}_{\!B}E \times_{B}E \cup \overset{k}{\product}_{B}E \times_{B} s(X)).
\end{align*}

Then the Whitehead definition of a fibrewise pointed L-S category $\catBB{E,s}$ for a fibrewise well-pointed space $E=(E,s)$ is given as follows:

\begin{defn}
$\catBB{E,s} \le m$ if the diagonal map $\Delta^{m+1} : E \to \overset{m+1}\product_{\!\!\!B} E$ is compressible into the fibrewise fat wedge by a fibrewise pointed homotopy.
\end{defn}

The fibrewise $A_{\infty}$-structure associated to the fibrewise loop space is fibrewise homotopy equivalent to the fibrewise version of the Ganea construction \cite{MR229240} as in \cite{MR2665235}, using the fibrewise version of fibre-cofibre lemma \cite{MR402694,MR1642747}.
More precisely, for a fibrewise well-pointed space $E=(E,s)$, the associated fibrewise projective $m$-space $P^{m}\LoopBB{E}$ together with a fibrewise fibration $E_{B}^{m+1}\LoopBB{E} \hookrightarrow P_{B}^{m}\LoopBB{E} \overset{e^{E}_{m}}\longrightarrow E$ and a fibrewise cofibration $E_{B}^{m+1}\LoopBB{E} \rightarrow P_{B}^{m}\LoopBB{E} \hookrightarrow P_{B}^{m+1}\LoopBB{E}$ is nothing but the $m$-th fibrewise Ganea construction, up to homotopy, and we obtain the following.

\begin{thm}\label{thm:fibrewisepointedsection}
$\catBB{E,s} \le m$ if and only if $e^{E}_{m} : P_{B}^{m}\LoopBB{E} \to E=(E,s)$ has a right homotopy inverse in the category of fibrewise pointed spaces and maps.
\end{thm}

We can define a fibrewise version of category and module weights as well, while they are more difficult to calculate than the usual category and module weights.
For a space $X$, we denote by $\double{X} = (\double{X},\proj_{1},X)$ the fibrewise space over $X$ as in \cite{MR2556074,MR2923451}.
If the space $X$ is good enough, we see that $s_{\id}=\Delta : X \hookrightarrow \double{X}$ is a closed cofibration, and hence $(\double{X},\Delta)$ is a fibrewise well-pointed space over $X$.

\begin{thm}[M.~Sakai and I]
$\tcm{X}=\catBB{\double{X},\Delta}$.
\end{thm}

The following result is recently announced in \cite{AX210911106}.

\begin{thm}[J. Aguilar-Guzmán and J. Gonz\'alez]
$\tcm{X}=\tc{X}$ if $X$ is ANR.
\end{thm}

%%%%%%%%%%%%%%%%%%%%%%%%%%%%%%%%%%%%%%%%%%%%%%%%%%%%%%%%%%%
\section{Fibrewise unpointed L-S theory}

Let us see TC from the fibrewise view point.
If there is a section of $\pi : \Path{X} \to \double{X}$ on a subset $F \subset \double{X}$, which gives a motion from the first point to the second point.
By reversing the orientation, we may think that a section gives a motion from the second point to the first point.
Thus by taking adjoint, we obtain a compression of $F$ into the diagonal $\Delta{X} \subset \double{X}$.
Now we assume that the inclusion $F \hookrightarrow \double{X}$ is a closed cofibration, and we obtain a deformation $h : \double{X} \!\times\! I \to \double{X}$ of the identity whose restriction to $F \!\times\! I$ coincides with the above compression of $F$ into $\Delta{X} \subset \double{X}$.
In particular, $h$ gives a motion of an element $(x,x) \in \Delta{X}$ to somewhere in $\Delta{X}$.

Let $E=(E,s)$ be a fibrewise pointed space over $X$.
If $E$ is covered by $m{+}1$ closed cofibrations each of which is equipped with a fibrewise compression into $s(X)$, then we obtain a fibrewise compression $H : E \times\! I \to \overset{m+1}\product_{\!\!\!B} E$ of the fibrewise diagonal $\Delta^{m+1}$ into the fibrewise fat wedge $\overset{m+1}{\fatvee}_{\!\!\!\!B}E$, where $H$ is a fibrewise unpointed homotopy.

\begin{defn}
$\catB{E,s} \le m$ if the fibrewise diagonal map $\Delta^{m+1} : E=(E,s) \to \overset{m+1}\product_{\!\!\!B}E$ is fibrewise compressible into the fibrewise fat wedge by a fibrewise unpointed homotopy.
\end{defn}

Then, M.~Sakai and the author obtained the following in \cite{MR2556074,MR2923451}.

\begin{thm}[M.~Sakai and I.]
$\catB{E,s} \le m$ if and only if $e^{E}_{m} : P_{B}^{m}\LoopBB{E} \to E=(E,s)$ has a right homotopy inverse in the category of fibrewise unpointed spaces and maps.
\end{thm}

For a space $X$, we denote by $\double{X} = (\double{X},\proj_{1},X)$ the fibrewise space over $X$, and $(\double{X},\Delta)$ is a fibrewise pointed space.
Then we obtain the following.

\begin{thm}[M.~Sakai and I.]\label{IS:tc-double}
$\tc{X}=\catB{\double{X},\Delta}$.
\end{thm}
It might be possible for us to extend Theorem \ref{IS:tc-double} to an equality of the $r$-th higher TC of a space $X$ and the fibrewise L-S category of a fibrewise space $(X^{r},\proj_{1},X)$ with the $r$-fold diagonal $\Delta_{r} : X \to X^{r}$ as its fibrewise base point.

Similarly to Theorem \ref{thm:fibrewisepointedsection}, we obtain the following.

\begin{thm}
$\tc{X} \le m$ if the natural map $e^{\double{X}}_{m} : P_{B}^{m}\LoopBB{\double{X}} \to \double{X} = (\double{X},\Delta)$ has a right homotopy inverse in the category of fibrewise unpointed spaces and maps.
\end{thm}

If we look at the fibrewise resolution of $\fLoop{BG}$ for the fibrewise L-S category, the cell-structure could be less complicated for a group $G$ than the fibrewise join for the sectional category.
For instance, for a Klein bottle $K$, D.~C.~Cohen and L.~Vandembroucq in \cite{MR3975552} obtained $\tc{K} = 4$ using a hard calculation of the obstruction produced from a fibrewise join (see A.~Costa and M.~Farber \cite{MR2649230}), while in \cite{MR3975098}, M.~Sakai, M.~Tsutaya and the author deduce the same result using a relatively shorter calculation. 
But the author has come to believe now that the fibrewise join of a projection is homotopy equivalent to the fibrewise projective space in some sense, and the fibrewise L-S theory and obstruction theory due to A.~Costa and M.~Farber based on a theorem of A.~S.~\v{S}varc are essentially the same: the difference is that the fibrewise L-S theory is more geometric and the obstruction theory is more group-theoretic.

\begin{prob}
Determine $\tc{M}$ for a topological spherical space form $M$.
\end{prob}
The answer to the above problem might suggest the way to determine the exact value of $\tc{\RP^{n}}=\Imm{\RP^{n}}$ for higher $n \ge 1$.

\appendix

%%%%%%%%%%%%%%%%%%%%%%%%%%%%%%%%%%%%%%%%%%%%%%%%%%%%%%%%%%%
\section{Filtrations of a classifying space}

Let $\Delta^{n}$ be the standard $n$-simplex.
$$
\Delta^{n} = \{(t_{0},\ldots,t_{n}) \mid t_{0}+\cdots+t_{n}=1\}.
$$
We can easily see that the above $\Delta^{n}$ is naturally homeomorphic to 
$$
{\Delta'}^{n} = \{(u_{0},\ldots,u_{n}) \mid 0 \le u_{0} \le \cdots \le u_{n}=1\}.
$$

%%%%%%
\subsection{Milnor construction}
We will explain briefly the construction introduced by J. Milnor \cite{MR77932}.
For a family of spaces $\mathbb{X} = \{X_{0}, \ldots, X_{n}\}$, we denote 
$$
\hat{E}_{M}(\mathbb{X}) = \Delta^{k} \times X_{0} \times \cdots \times X_{n}/\sim
$$
where the equivalence relation $\sim$ is generated by 
\begin{align*}&
(t_{0},\ldots,t_{i-1},0,t_{i+1}\ldots,t_{n};x_{0},\ldots,x_{i-1},x_{i},x_{i+1},\ldots,x_{n})
\\[0ex]&\hskip5em\sim (t_{0},\ldots,t_{i-1},0,t_{i+1}\ldots,t_{n};x_{1},\ldots,x_{i-1},x'_{i},x_{i+1},\ldots,x_{n}),\quad 0 \!\le\! i \!\le\! n.
\end{align*}\vskip-.5ex\noindent
The class of $(t_{0},\ldots,t_{k};x_{0},\ldots,x_{k})$ is denoted by $\underset{i=0}{\overset{k}\sum}\,t_{i}{\cdot}x_{i}$, in which $0{\cdot}x_{i}$ is often ignored.
So, we may assume that $\hat{E}_{M}(\mathbb{X}) \subset \hat{E}_{M}(\mathbb{X}')$ if $\mathbb{X}$ is a sub-family of $\mathbb{X}'$.

When $\mathbb{X}$ is given by $X_{i}=X$, $0 \le i \le n$, for a fixed space $X$, we denote $\hat{E}_{M}^{\,n}(X)=\hat{E}_{\!M}(\mathbb{X})$, in this section.
If further $X$ has a homotopy type of a CW complex, then $\hat{E}_{\!M}^{\infty}(X) = \underset{n}\bigcup \hat{E}_{M}^{\,n+1}(X)$ is contractible (see J.~Milnor \cite{MR100267}).

For a topological group $G$, $\hat{E}_{M}^{\,n}(G)$ and $\hat{E}_{\!M}^{\infty}(G)$ are free $G$ spaces by the diagonal $G$ action.
Thus we obtain a fibration $G \hookrightarrow \hat{E}_{M}^{\,n}(G) \twoheadrightarrow \hat{P}_{M}^{\,n-1}(G)$ over $\hat{P}_{M}^{\,n-1}(G)=\hat{E}_{M}^{\,n}(G)/G$ which is a subspace of $\hat{P}_{M}^{\infty}(G)=\hat{E}_{M}^{\infty}(G)/G$.
The fibration is induced from the universal fibration $\hat{E}_{M}^{\infty}(G) \twoheadrightarrow \hat{P}_{M}^{\infty}(G)$, and we obtain the following fibre sequences.
\begin{equation}
\hat{E}_{M}^{\,n}(G) \twoheadrightarrow \hat{P}_{M}^{\,n-1}(G) \to \hat{P}_{M}^{\infty}(G),\quad n \ge 1.\tag{F}
\end{equation}
Moreover by definition, we have the following \textit{unreduced} cofibre sequences.
\begin{align}&
\hat{P}_{M}^{\,0}(G) = \{\ast\},\tag{Z}
\\&
\hat{E}_{M}^{\,n}(G) \twoheadrightarrow \hat{P}_{M}^{\,n-1}(G) \hookrightarrow \hat{P}_{M}^{\,n}(G),\quad n \ge 1.\tag{C}
\end{align}

\begin{expl}
For $S^{0}\!=\!\mathrm{O}(1)$, $S^{1}\!=\!\mathrm{U}(1)$ and $S^{3}\!=\!\mathrm{Sp}(1)$, the following holds.
\begin{enumerate}
\item
$\hat{P}_{M}^{\,n}(S^{0}) \homeo \RP^{n}$ the real projective $n$-space.\vspace{.5ex}
\item
$\hat{P}_{M}^{\,n}(S^{1}) \homeo \CP^{n}$ the complex projective $n$-space.\vspace{.5ex}
\item
$\hat{P}_{M}^{\,n}(S^{3}) \homeo \HP^{n}$ the quaternionic projective $n$-space.
\end{enumerate}
\end{expl}

%%%%%%
\subsection{Geometric Bar construction}
Let $G$ be a topological group with unit $e$ as the non-degenerate base point.
For a right $G$-space $X$ and a left $G$-space $Y$, we have the bar construction $E_{M}^{\,n+1}(G)$ and $P_{M}^{\,n}(G)$ which are \textit{reduced} versions of Milnor's construction $\hat{E}^{\,n+1}_{M}(G)$ and $\hat{P}^{\,n}_{M}(G)$ (see R.~J.~Milgram \cite{MR208595}):
\begin{align*}&
E_{M}^{\,n+1}(G) = \underset{0 \le k \le n}{\textstyle\coprod}\,\Delta^{k} \times (G \times G^{k})/\sim',
\qquad%
P_{M}^{\,n}(G) = \underset{0 \le k \le n}{\textstyle\coprod}\,\Delta^{k} \times (\ast \times\ G^{k})/\sim'
\\[.5ex]&\quad
(t_{0},\ldots,t_{i-1},0,t_{i+1}\ldots,t_{k};x,g_{1}\ldots,g_{i},\ldots,g_{k})
\\[.25ex]%
&\hskip5em\sim' \begin{cases}
(t_{0},\ldots,t_{k-1};xg_{1},g_{2},\ldots,g_{k}),&i \!=\! 0,
\\[.25ex]%
(t_{0},\ldots,t_{i-1},t_{i+1}\ldots,t_{k};x,g_{1},\ldots,g_{i}g_{i+1},\ldots,g_{k}dy),&0 \!<\! i \!<\! k,
\\[.25ex]%
(t_{0},\ldots,t_{k-1};x,g_{1},\ldots,g_{k-1}),&i \!=\! k,
\end{cases}
\\[.5ex]&\quad
(t_{0},\ldots,t_{i-1},t_{i},t_{i+1},\ldots,t_{k};x,g_{1},\ldots,g_{i-1},e,g_{i+1},\ldots,g_{k})\quad (1 \!\le\! i \!\le\! k)
\\[0ex]&\hskip5em\sim' (t_{0},\ldots,t_{i-1}{+}t_{i},t_{i+1}\ldots,t_{k};x,g_{1},\ldots,g_{i-1},g_{i+1},\ldots,g_{k}).
\end{align*}\vskip.5ex\noindent
where $xg$ is the group multiplication if $X=G$, and ${\ast}g=\ast$ if $X=\ast$.

Then we have the following relative homeomorphisms.
\begin{align*}&
(\Delta^{n} \times G{\times}G^{n}, \partial\Delta^{n} \times  G \times G^{n} \cup \Delta^{n} \times G \times G^{[n]}) \twoheadrightarrow (E_{M}^{\,n+1}(G),E_{M}^{\,n}(G)),
\\&
(\Delta^{n} \times G^{n}, \partial\Delta^{n} \times  G^{n} \cup \Delta^{n} \times G^{[n]}) \twoheadrightarrow (P_{M}^{\,n}(G),P_{M}^{\,n-1}(G)),
\end{align*}
where $G^{[n]} = \{ (g_{1},\dots,g_{n}) \mid \exists\,i \ g_{i} \!=\! e \}$.
The free action of $G$ on $E_{M}^{\,n}(G)$ is given by
$$
h{\cdot}(t,t_{1},\ldots,t_{k};g,g_{1},\ldots,g_{k},\ast) = (t,t_{1},\ldots,t_{k};hg,g_{1},\ldots,g_{k},\ast).
$$
If further $G$ has a homotopy type of a CW complex, then ${E}_{\!M}^{\infty}(G) = \underset{n}\bigcup\,{E}_{M}^{\,n+1}(G)$ is contractible (see J.~Milnor \cite{MR100267}).
Thus we obtain a fibration $G \hookrightarrow {E}_{M}^{\,n}(G) \twoheadrightarrow {P}_{M}^{\,n-1}(G)$ over ${P}_{M}^{\,n-1}(G) \subset {P}_{M}^{\infty}(G)$, which is induced from the universal fibration ${E}_{M}^{\infty}(G) \twoheadrightarrow {P}_{M}^{\infty}(G)$, so that we obtain the following fibre sequences.
\begin{equation}
{E}_{M}^{\,n}(G) \twoheadrightarrow {P}_{M}^{\,n-1}(G) \to {P}_{M}^{\infty}(G),\quad n \ge 1.\tag{F}
\end{equation}
Moreover by definition, we have the following cofibre sequences.
\begin{align}&
{P}_{M}^{\,0}(G) = \{\ast\},\tag{Z}
\\&
{E}_{M}^{\,n}(G) \twoheadrightarrow {P}_{M}^{\,n-1}(G) \hookrightarrow {P}_{M}^{\,n}(G),\quad n \ge 1.\tag{C}
\end{align}

\begin{rem}
$P_{M}^{\,\infty}(G)$ and $\hat{P}_{M}^{\,\infty}(G)$ have the same homotopy type.
Moreover, for all $n \!\ge\! 1$, $P_{M}^{\,n}(G)$ and $\hat{P}_{M}^{\,n}(G)$ have the same homotopy type as well.
\end{rem}

%%%%%%%%%%%%%%%%%%%%%%%%%%%%%%%%%%%%%%%%%%%%%%%%%%%%%%%%%%%
\section{$A_{\infty}$-forms and $A_{\infty}$-structures}\label{app-sect:higherassociativity}

\subsection{$A_{\infty}$-form and $A_{\infty}$-structure for a space}

First, we denote the associahedron by $K(n)$, $n \ge 2$, which is given by
\begin{align*}&
K(n)=\{ (t_{1},\dots,t_{n}) \in \real_{+}^{n}\mid \forall\,i \ 0 \le t_{1}\!+\cdots+\!t_{i} \le i{-}1, \ t_{1}\!+\cdots+\!t_{n}=n{-}1 \},
\end{align*}
where $\real_{+}=[0,\infty)$, which is slightly modified from the original construction by J. D. Stasheff \cite{MR0158400}.
Clearly, the above $K(n)$ can naturally be identified with 
$$
K'(n) = \{(u_{1},\ldots,u_{n}) \in \real^{n} \mid 0 \!=\! u_{1} \!\le\! u_{2} \!\le\! \cdots \!\le\! u_{n-1} \!\le\! u_{n} \!=\! n{-}1, \ \forall\,i \ \,u_{i} \le i{-}1\}.
$$
Let $A(n)=\{(k,r,s) \mid 1 \!\le\! k \!\le\! r, \ 2 \!\le\! s\!=\!n{+}1{-}r \!\le\! n{-}1\}$.
Then we have maps
\begin{align*}&
\partial_{k} : K(r) \times K(s) \to \partial K(n) \subset K(n), \hskip-3em &&(k,r,s) \in A(n)\quad\text{and}
\\&
s_{j} : K(n) \to K(n{-}1), \hskip-3em &&1 \!\le\! j \!\le\! n,
\end{align*}
which are defined as follows (see Stasheff \cite{MR0158400} or \cite{AX12115741}).
\begin{align*}&
\partial_{k}((t_{1},\ldots,t_{r}),(u_{1},\ldots,u_{s})) = 
(t_{1},\ldots,t_{k-1},u_{1},\ldots,u_{s}{+}t_{k},\ldots,t_{r}),\quad 1\!\le\!k\!\le\!r,
\\&
s_{j}(t_{1},\ldots,t_{n})=
\left\{\vphantom{\Bigg(}\!\!\!\!\right.\begin{array}{ll}
(t'_{2},\ldots,t'_{n}), \ \ \xi(t_{1},\ldots,t_{n})\!=\!(0,t'_{2},\ldots,t'_{n}), & j\!=\!1,
\\[.5ex]
(t_{1},\ldots,t_{j-1}{+}t'_{j},\ldots,t'_{n}), \ \ \xi(t_{j},\ldots,t_{n})\!=\!(t'_{j},\ldots,t'_{n}), & 1\!<\!j\!\le\!n,
\end{array}
\end{align*} 
where $\xi : \real_{+}^{n} \ni (t_{1},\ldots,t_{n}) \mapsto (t'_{1},\ldots,t'_{n}) \in \real_{+}^{n}$ is given by %
\begin{align*}&
t'_{k}=\begin{cases}
\max\{0,t_{1}{-}1\}, & k\!=\!1,
\\[.5ex]
\min\left\{\vphantom{\bigg(}\!\right. t_{k}, \underset{1 \le j \le k}\max \left\{\sum_{i=1}^{j}(t_{i}{-}1)\right\} - \sum_{i=1}^{k-1}(t'_{i}{-}1)\left.\!\vphantom{\bigg)}\right\}, & 1\!<\!k\!\le\!n.
\end{cases}
\end{align*}\vskip1ex\noindent
By definition, they satisfy the following formula for $(k,r,s) \!\in\! A(n)$ and $1 \!\le\! j \!\le\! n$. 
$$
s_{j}\comp\partial_{k}(\rho,\sigma) = \begin{cases}
\partial_{k-1}{\comp}(s_{j}(\rho) \times \sigma),& j \!<\! k, \ r\!>\!2,
\\[-.5ex]
\sigma,& j\!=\!1, \,k\!=\!2, \,r\!=\!2,
\\
\partial_{k}{\comp}(\rho \times s_{j-k+1}(\sigma)),& k \!\le\! j \!<\! k{+}s, \ r \!<\! n{-}1,
\\[-.5ex]
\rho,& k\!\le\!j\!\le\!k{+}1, \ r\!=\!n{-}1,
\\
\partial_{k-1}{\comp}(s_{j-s+1}(\rho) \times \sigma),& k{+}s \!\le\! j \!\le\! n, \ r\!>\!2,
\\[-.5ex]
\sigma,& j\!=n, \,k\!=\!1, \,r\!=\!2.
\end{cases}
$$

\smallskip

An $A_{\infty}$-form on a space $G$ is a sequence of maps
$$
\alpha_{n} : K(n) \times G^{n} \to G
$$
satisfying boundary and unital conditions for $(k,r,s) \!\in\! A(n)$ and $1 \!\le\! j \!\le\! n$:
\begin{enumerate}
\item
$\alpha_{n}(\partial_{k}(\rho,\sigma);g_{1},\ldots,g_{n})=\alpha_{r}(\rho;g_{1},\ldots,\alpha_{s}(\sigma;g_{k},\ldots,g_{k+s-1}),\ldots,g_{n})$.\vspace{1ex}
\item
$\alpha_{n}(\tau;g_{1},\ldots,g_{j-1},e,g_{j+1},\ldots,g_{n})
\\\hphantom{\qquad\qquad\qquad\qquad} =\alpha_{n-1}(s_{j}(\tau);g_{1},\ldots,g_{j-1},g_{j+1},\ldots,g_{n})$.
\end{enumerate}

\smallskip

A space with $A_{\infty}$-form is called an $A_{\infty}$-space.
For an $A_{\infty}$-space $G$, we have $E^{\,n+1}(G)$ and $P^{\,n}(G)$ the standard $A_{\infty}$-structure of $G$, which are inductively defined by the following relative homeomorphisms
\begin{align*}&
G \!\times\! (K(n{+}2) \!\times\! G^{n}, \partial K(n{+}2) \!\times\! G^{n} \!\cup\! K(n{+}2) \!\times\! G^{[n]}) \overset{\phi_{n}}\twoheadlongrightarrow (E^{n+1}(G),E^{n}(G)) \ \ \text{and}
\\&
\ast \!\times\!\ (K(n{+}2) \!\times\! G^{n}, \partial K(n{+}2) \!\times\! G^{n} \!\cup\! K(n{+}2) \!\times\! G^{[n]}) \overset{\phi_{n}}\twoheadlongrightarrow (P^{n}(G),P^{n-1}(G))
\end{align*}
which is defined on $X \!\times\! (\partial K(n{+}2) \!\times\! G^{n} \!\cup\! K(n{+}2) \!\times\! G^{[n]})$, $X\!=\!G$ or $\ast$, as follows (see Stasheff \cite{MR0158400}, \cite{MR1000378} or \cite{AX12115741}):
for $\mathbold{x}\!=\!(x;\partial_{k}(\rho,\sigma);g_{2},\ldots,g_{n+1}) \!\in\! X \!\times\! \partial K(n{+}2) \!\times\! G^{n}$, $(k,r,s) \!\in\! A(n)$, 
\begin{align*}&
\phi_{n}(\mathbold{x})=\begin{cases}\,
\phi_{r}(\alpha'_{s}(\sigma;x,g_{2},\ldots,g_{s});\rho;g_{s+1},\ldots,g_{n+1}),&k\!=\!1,
\\[.5ex]\,
\phi_{r}(x;\rho;g_{1},\ldots,\alpha_{s}(\sigma;g_{k},\ldots,g_{k+s-1}),\ldots,g_{n+1}),&1\!<\!k\!<\!r,
\\[.5ex]\,
\phi_{r}(x;\rho;g_{1},\ldots,g_{r-1}),&k\!=\!r,
\end{cases}\quad \text{and}
\intertext{for $\mathbold{y}=(x;\tau;g_{2},\ldots,g_{j-1},e,g_{j+1},\ldots,g_{n+1}) \in X \times K(n{+}2) \times G^{[n]}$, $1 \!<\! j \!<\! n{+}2$,}&
\phi_{n}(\mathbold{y}) = \phi_{n-1}(x;s_{j}(\tau);g_{2},\ldots,g_{j-1},g_{j+1},\ldots,g_{n+1}),
\end{align*}
where $\alpha'_{s}=\alpha_{s}$ if $X=G$, and $\alpha'_{s}(\sigma;\ast,g_{2},\ldots,g_{s})$ $=$ $\ast$ if $X$ $=$ $\ast$.

\smallskip

Then the canonical projection $G \!\times\! K(n{+}2) \!\times\! G^{n}$ $\to$ $\{\ast\} \!\times\! K(n{+}2) \!\times\! G^{n}$ induces a projection $p_{n}^{G} : E^{n+1}(G) \to P^{n}(G)$.
We also have a contractible subspace $D^{n}(G)$ of $E^{n+1}(G)$ obtained by the relative homeomorphism
\begin{align*}&
(K(n{+}2) \!\times\! G^{n}, \partial K(n{+}2) \!\times\! G^{n} \!\cup\! K(n{+}2) \!\times\! G^{[n]}) \overset{\psi_{n}}\twoheadlongrightarrow (D^{n}(G),E^{n}(G)) 
\end{align*}
which is defined on $\partial K(n{+}2) \!\times\! G^{n} \!\cup\! K(n{+}2) \!\times\! G^{[n]}$ as follows (see J.~D.~Stasheff \cite{MR0158400}, \cite{MR1000378} or \cite{AX12115741}):
for $\mathbold{x} \!=\! (\partial_{k}(\rho,\sigma);g_{2},\ldots,g_{n+1}) \!\in\! \partial K(n{+}2) \!\times\! G^{n}$, $(k,r,s) \!\in\! A(n)$, 
\begin{align*}&\,
\psi_{n}(\mathbold{x})=\begin{cases}\,
\phi_{r}(\alpha_{s}(\sigma;e,g_{2},\ldots,g_{s});\rho;g_{s+1},\ldots,g_{n+1}),&k\!=\!1,
\\[.5ex]\,
\phi_{r}(e;\rho;g_{2},\ldots,\alpha_{s}(\sigma;g_{k},\ldots,g_{k+s-1}),\ldots,g_{n+1}),&1\!<\!k\!<\!r,
\\[.5ex]\,
\phi_{r}(e;\rho;g_{2},\ldots,g_{r-1}),&k\!=\!r,
\end{cases}\quad \text{and}
\intertext{for $\mathbold{y}=(\tau;g_{2},\ldots,g_{j-1},e,g_{j+1},\ldots,g_{n+1}) \in K(n{+}2) \times G^{[n]}$, $1\!<\!j\!<\!n{+}2$,}&\,
\psi_{n}(\mathbold{y})=\phi_{n-1}(e;s_{j}(\tau);g_{2},\ldots,g_{j-1},g_{j+1},\ldots,g_{n+1}).
\end{align*}

\smallskip

We obtain the following cofibre sequences, since $D^{n}(G)$ is contractible:
\begin{align}&
P^{0}(G) = \{\ast\},\tag{Z}
\\&
E^{n}(G) \twoheadrightarrow P^{n-1}(G) \hookrightarrow P^{n}(G),\quad n \ge 1.\tag{C}
\end{align}
Since the fibration $p^{G}_{n} : E^{n}(G) \twoheadrightarrow P^{n-1}(G)$ is a pull-back of the universal fibration $E^{\infty}(G) \to P^{\infty}(G)$ with $E^{\infty}(G) \simeq \ast$, we also obtain the following fibre sequences:
\begin{equation}
E^{n}(G) \twoheadrightarrow P^{n-1}(G) \hookrightarrow P^{\infty}(G),\quad n \ge 1,\tag{F}
\end{equation}
where the fibre of $p^{G}_{n} : E^{n}(G) \twoheadrightarrow P^{n-1}(G)$ is $G$.
\begin{rem}
The projective $n$-space $P^{n}(G)$ of an $A_{\infty}$-space $G$ is denoted in \cite{MR0158400} by $GP^{n}$ and is called the $G$-projective $n$-space.
For example, $\real P^{n} \simeq S^{0}P^{n}=P^{n}(S^{0})$ and $\complex P^{n} \simeq S^{1}P^{n}=P^{n}(S^{1})$, and etc.
\end{rem}

\subsection{$A_{\infty}$-form and $A_{\infty}$-structure for a map}

Second, we denote the multiplihedron by $J(n)$, $n \ge 1$, which is given by
\begin{align*}&
J(n)=\{ (t_{1},\dots,t_{n}) \in \real_{+}^{n}\mid \forall\,i \ 0\!\le\!t_{1}\!+\cdots+\!t_{i}\!\le\!i{-}1{+}a, \ t_{1}\!+\cdots+\!t_{n} \!=\! n{-}1{+}a \}, 
\end{align*}
where $a\!=\!\fracinline1/2$.
Clearly, the above $J(n)$ can naturally be identified with 
$$
J'(n) = \{(u_{1},\ldots,u_{n}) \in \real^{n} \mid 0 \!\le\! u_{1} \!\le\! \cdots \!\le\! u_{n} \!=\! n{-}1{+}a, \ \forall\,i \ 0 \!\le\! u_{i} \!\le\! i{-}1{+}a\}.
$$
Let $A'(n)=\{(k,r,s) \mid 1 \!\le\! k \!\le\! r, \ 2 \!\le\! s\!=\!n{+}1{-}r \!\le\! n\}$, $B(n)=\{(t;r_{1},\dots,r_{t})$ $\mid$ $2$ \!$\le$\! $t \!\le\! n, \ \forall\,i \ 1 \!\le\! r_{i} \!<\! n, \ r_{1}{+}\cdots{+}r_{t}\!=\!n\}$. 
Then we have maps
\begin{align*}&
\delta_{k} : J(r) \times K(s) \to \partial J(n) \subset J(n), \ (k,r,s) \in A'(n),
\\&
\delta : K(t) \times J(r_{1}) \times \cdots \times J(r_{t}) \to \partial J(n) \subset J(n), \ (t;r_{1},\dots,r_{t}) \in B(n)\quad\text{and}
\\&
d_{j} : J(n) \to J(n{-}1), \ 1 \!\le\! j \!\le\! n,
\end{align*}
which are defined as follows (see \cite{MR1000378} or  \cite{AX12115741}).
\begin{align*}&
\delta_{k}(v_{1},\ldots,v_{r};u_{1},\ldots,u_{s})=(v_{1},\dots,v_{k-1},u_{1},\dots,u_{s}{+}v_{k},\dots,v_{r}),\quad 1 \!\le\! k \!\le\! r,
\\[.5ex]&
\delta(u_{1},\ldots,u_{t};v^{(1)}_{1}\!,\ldots,v^{(1)}_{r_{1}},\ldots,v^{(t)}_{1},\ldots,v^{(t)}_{r_{t}})
\\&\hphantom{\qquad}
=(v^{(1)}_{1},\ldots,v^{(1)}_{r_{1}-1},v^{(1)}_{r_{1}}{+}(1{-}a){\cdot}u_{1},\ldots,v^{(t)}_{1},\ldots,v^{(t)}_{r_{t}-1},v^{(t)}_{r_{t}}{+}(1{-}a){\cdot}u_{t}),
\\[.5ex]&
d_{j}(t_{1},\ldots,t_{n})=
\!\left\{\vphantom{\Bigg(}\!\!\!\!\!\right.\begin{array}{ll}
(t'_{2},\ldots,t'_{n}), \ \ \xi(t_{1},\ldots,t_{n})\!=\!(0,t'_{2},\ldots,t'_{n}), & \!j\!=\!1,
\\[.5ex]
(t_{1},\ldots,t_{j-1}{+}t'_{j},\ldots,t'_{n}), \ \ \xi(t_{j},\ldots,t_{n})\!=\!(t'_{j},\ldots,t'_{n}), & \!1\!<\!j\!\le\!n.
\end{array}
\end{align*}
By definition, they satisfy the following formulas for $(k,r,s) \!\in\! A'(n)$, $(t;r_{1},\ldots,r_{t})$ \!$\in$\! $B(n)$ and $1 \!\le\! j \!\le\! n$.
\begin{align*}&
d_{j}\comp\delta_{k}(\rho,\sigma) = \begin{cases}
\delta_{k-1}{\comp}(d_{j}(\rho) \times \sigma),& j \!<\! k,
\\[-.25ex]
\delta_{k}{\comp}(\rho \times s_{j-k+1}(\sigma)),& k \!\le\! j \!<\! k{+}s, \ r \!<\! n{-}1,
\\[-.5ex]
\rho,& k\!\le\!j\!\le\!k{+}1, \ r\!=\!n{-}1,
\\[-.25ex]
\delta_{k-1}{\comp}(d_{j-s+1}(\rho) \times \sigma),& k{+}s \!\le\! j \!\le\! n,
\end{cases}
\\&
d_{j}{\comp}\delta(\tau;\rho_{1},\ldots,\rho_{t})= \begin{cases}
\delta(\tau;\rho_{1},{\cdots},
d_{j-\hat{r}_{k-1}}(\rho_{k}),{\cdots},\rho_{t}), & \hat{r}_{k-1} \!<\! j \!\leq\! \hat{r}_{k}, \ r_{k}\!>\!1,
\\
\delta(s_{k}(\tau);\rho_{1},{\cdots},\rho_{k-1},\rho_{k+1},{\cdots},\rho_{t}), & j\!=\!\hat{r}_{k}, \ r_{k}\!=\!1, \ t\!>\!2,
\\[-.5ex]
\rho_{2}, & j\!=\!1, \ r_{1}\!=\!1, \ t\!=\!2,
\\[-.5ex]
\rho_{1}, & j\!=\!n, \ r_{2}\!=\!1, \ t\!=\!2,
\end{cases}
\end{align*}
where $\hat{r}_{0}=0$ and $\hat{r}_{k}=r_{1}+\cdots+r_{k}$, $1 \!\le\! k \!\le\! t$.

\smallskip

An $A_{\infty}$-form on a map $f : G \to H$ of $A_{\infty}$-spaces is a sequence of maps
$$
\beta_{n} : J(n) \times G^{n} \to H
$$
which satisfies the boundary and the unital conditions for $(k,r,s) \!\in\! A'(n)$, $(t;r_{1},\ldots,r_{t})$ \!$\in$\! $B(n)$ and $1 \!\le\! j \!\le\! n$ as follows (see \cite{MR1000378} or \cite{AX12115741}):
\begin{enumerate}
\item
$\beta_{n}(\delta_{k}(\rho,\sigma);g_{1},\ldots,g_{n})=\beta_{r}(\rho;g_{1},\ldots,\alpha^{G}_{s}(\sigma;g_{k},\ldots,g_{k+s-1}),\ldots,g_{n})$.\vspace{.5ex}
\item
$\beta_{n}(\delta(\tau;\rho_{1},\ldots,\rho_{t});g_{1},\ldots,g_{n})$ \vspace{.5ex}\\\hphantom{==} $=$ $\alpha^{H}_{r}(\tau;\beta_{r_{1}}(\rho_{1};g_{1},\ldots,g_{r_{1}}),\ldots,\beta_{r_{t}}(\rho_{t};g_{n-r_{t}+1},\ldots,g_{n}))$.\vspace{.5ex}
\item
$\beta_{n}(\tau;g_{1},\ldots,g_{j-1},e,g_{j+1},\ldots,g_{n})$ \vspace{.5ex}\\\hphantom{\qquad}$=$ $\beta_{n-1}(d_{j}(\tau);g_{1},\ldots,g_{j-1},g_{j+1},\ldots,g_{n})$, \ \ $1 \!\le\! j \!<\! n$,
\end{enumerate}
where $\{\alpha^{G}_{n}\}$ and $\{\alpha^{H}_{n}\}$ denote the $A_{\infty}$-structures on $G$ and $H$, respectively.

\smallskip

A map with an $A_{\infty}$-form is called an $A_{\infty}$-map.
For an $A_{\infty}$-map $f$, the associated $A_{\infty}$-structure $(E^{n}(f),P^{n-1}(f))$ is a pair of maps given to fit in with the following natural commutative diagram for $n \!\ge\! 1$ (see \cite{MR1000378} or \cite{AX12115741}).

$$
\begin{tikzcd}
{E^{n}(G)}
\arrow[rrr,"E^{n}(f)"]
\arrow[ddd,"p^{G}_{n}"']
\arrow[rd,hook]
&{}
&{}
&{E^{n}(H)}
\arrow[ddd,"p^{H}_{n}"]
\arrow[rd,hook]
\\
{}
&{D^{n}(G)}
\arrow[rrr,"D^{n}(f)"]
\arrow[ddd,"q^{G}_{n}"']
&{}
&{}
&{D^{n}(H)}
\arrow[ddd,"q^{H}_{n}"']
\\
\\
{P^{n-1}(G)}
\arrow[rrr,"P^{n-1}(f)"']
\arrow[rd,hook]
&{}
&{}
&{P^{n-1}(H)}
\arrow[rd,hook]
\\
{}
&{P^{n}(G)}
\arrow[rrr,"P^{n}(f)"']
&{}
&{}
&{P^{n}(H).\!}
\end{tikzcd}
$$

%%%%%%%%%%%%%%%%%%%%%%%%%%%%%%%%%%%%%%%%%%%%%%%%%%%%%%%%%%%
\section{Ganea's fibre-cofibre construction}\label{app-sect:ganea-fibre-cofibre}

T.~Ganea introduced his (fibre-cofibre) construction to show the characterization of the L-S category \cite[Proposition 2.2]{MR229240}.
Ganea began with (Z) below and inductively constructed a series of cofibrations (C) and fibrations (F) concretely.
\begin{align}&
E_{0}(X) = \{\ast\} \overset{p_{0}}\longrightarrow X,
\tag{Z}
\\&
\text{$F_{n}(X)$ is the mapping fibre of $p_{n} : E_{n}(X) \to X$,}
\tag{F}
\\&
\text{$E_{n+1}(X)$ is the mapping cone of $i_{n} : F_{n}(X) \hookrightarrow E_{n}(X)$.}
\tag{C}
\end{align}
Here, $p_{n+1}$ is the extension of $p_{n} : E_{n}(X) \to X$ by annihilating whole $C(F_{n}(X))$, which is well-defined since $p_{n}(F_{n}(X))=\{\ast\}$.
This construction is, however, depending on the choice of $p_{n+1}$, since we know there are many choices for $p_{n+1}$ extending $p_{n}$. In any case, Exercise \ref{exer:ganeaspaces} tells us that such choices does not affect the homotopy type of fibre-cofibre constructions. But this is not so obvious as we have seen in the fact that \cite[Exercise 2.16]{MR1990857} can not be verified as is.

So, it would be possible to say that the space $E_{n}(X)$, which is often referred as $G_{n}(X)$ the Ganea space, is homotopy equivalent to $\Loop{X}$-projective $n$-space $P^{n}\Loop{X}$, which is obtained by using Exercise \ref{exer:ganeaspaces} and Theorem \ref{thm:stasheff-main}.

%%%%%%%%%%%%%%%%%%%%%%%%%%%%%%%%%%%%%%%%%%%%%%%%%%%%%%%%%%%
\section*{Acknowledgements}
The author would like to thank the editors for allowing him to discuss the works of E.~Fadell and S.~Husseini from the $A_{\infty}$-theoretic view point, and the referees for their valuable comments on earlier versions of this article.
This research was supported by Grant-in-Aid for Scientific Research (S) \#17H06128 and for Exploratory Research \#18K18713 from Japan Society for the Promotion of Science.

%%%%%%%%%%%%%%%

\bibliographystyle{alpha}
\bibliography{lscat2}

\begin{thebibliography}{FSGTST04}

\bibitem[AGG21]{AX210911106}
Jorge Aguilar-Guzm{\'a}n and Jes{\'u}s Gonz{\'a}lez.
\newblock Motion planning in polyhedral products of groups and a
  fadell-husseini approach to topological complexity, 2021.
\newblock arXiv: https://arxiv.org/abs/2109.11106.

\bibitem[Ber58]{MR96221}
Isra\"{e}l Ber\v{s}te\u{\i}n.
\newblock Sur la cat\'{e}gorie de {L}usternik-{S}chnirelmann.
\newblock {\em C. R. Acad. Sci. Paris}, 246:362--364, 1958.

\bibitem[Ber76]{MR400212}
Israel Berstein.
\newblock On the {L}usternik-{S}chnirelmann category of {G}rassmannians.
\newblock {\em Math. Proc. Cambridge Philos. Soc.}, 79(1):129--134, 1976.

\bibitem[BG62]{MR139168}
I.~Berstein and T.~Ganea.
\newblock The category of a map and of a cohomology class.
\newblock {\em Fund. Math.}, 50:265--279, 1961/62.

\bibitem[BH60]{MR126276}
I.~Berstein and P.~J. Hilton.
\newblock Category and generalized {H}opf invariants.
\newblock {\em Illinois J. Math.}, 4:437--451, 1960.

\bibitem[CF10]{MR2649230}
Armindo Costa and Michael Farber.
\newblock Motion planning in spaces with small fundamental groups.
\newblock {\em Commun. Contemp. Math.}, 12(1):107--119, 2010.

\bibitem[CJ98]{MR1646248}
Michael Crabb and Ioan James.
\newblock {\em Fibrewise homotopy theory}.
\newblock Springer Monographs in Mathematics. Springer-Verlag London, Ltd.,
  London, 1998.

\bibitem[CLOT03]{MR1990857}
Octav Cornea, Gregory Lupton, John Oprea, and Daniel Tanr\'{e}.
\newblock {\em Lusternik-{S}chnirelmann category}, volume 103 of {\em
  Mathematical Surveys and Monographs}.
\newblock American Mathematical Society, Providence, RI, 2003.

\bibitem[Cor94]{MR1259517}
Octavian Cornea.
\newblock Cone-length and {L}usternik-{S}chnirelmann category.
\newblock {\em Topology}, 33(1):95--111, 1994.

\bibitem[CV17]{MR3975552}
Daniel~C. Cohen and Lucile Vandembroucq.
\newblock Topological complexity of the {K}lein bottle.
\newblock {\em J. Appl. Comput. Topol.}, 1(2):199--213, 2017.

\bibitem[DR09]{MR2475974}
Alexander~N. Dranishnikov and Yuli~B. Rudyak.
\newblock On the {B}erstein-{S}varc theorem in dimension 2.
\newblock {\em Math. Proc. Cambridge Philos. Soc.}, 146(2):407--413, 2009.

\bibitem[EG57]{MR85510}
Samuel Eilenberg and Tudor Ganea.
\newblock On the {L}usternik-{S}chnirelmann category of abstract groups.
\newblock {\em Ann. of Math. (2)}, 65:517--518, 1957.

\bibitem[Far03]{MR1957228}
Michael Farber.
\newblock Topological complexity of motion planning.
\newblock {\em Discrete Comput. Geom.}, 29(2):211--221, 2003.

\bibitem[FG08]{MR2407101}
Michael Farber and Mark Grant.
\newblock Robot motion planning, weights of cohomology classes, and cohomology
  operations.
\newblock {\em Proc. Amer. Math. Soc.}, 136(9):3339--3349, 2008.

\bibitem[FH92]{MR1317569}
Edward Fadell and Sufian Husseini.
\newblock Category weight and {S}teenrod operations.
\newblock volume~37, pages 151--161. 1992.
\newblock Papers in honor of Jos\'{e} Adem (Spanish).

\bibitem[FH94]{MR1397466}
Edward Fadell and Sufian~Y. Husseini.
\newblock Relative category, products and coproducts.
\newblock {\em Rend. Sem. Mat. Fis. Milano}, 64:99--115 (1996), 1994.

\bibitem[FO99]{MR1688434}
Kenji Fukaya and Kaoru Ono.
\newblock Arnold conjecture and {G}romov-{W}itten invariant.
\newblock {\em Topology}, 38(5):933--1048, 1999.

\bibitem[Fox39]{MR2937306}
Ralph~H. Fox.
\newblock {\em O{N} {THE} {LUSTERNIK} {SCHNIRELMANN} {CATEGORY}}.
\newblock PhD thesis, 1939.
\newblock Thesis (Ph.D.)--Princeton University.

\bibitem[Fox41]{MR4108}
Ralph~H. Fox.
\newblock On the {L}usternik-{S}chnirelmann category.
\newblock {\em Ann. of Math. (2)}, 42:333--370, 1941.

\bibitem[FSGTST04]{MR2022385}
Luc\'{\i}a Fern\'{a}ndez-Su\'{a}rez, Antonio G\'{o}mez-Tato, Jeffrey Strom, and
  Daniel Tanr\'{e}.
\newblock The {L}usternik-{S}chnirelmann category of {$\rm Sp(3)$}.
\newblock {\em Proc. Amer. Math. Soc.}, 132(2):587--595, 2004.

\bibitem[Gan67]{MR229240}
T.~Ganea.
\newblock Lusternik-{S}chnirelmann category and strong category.
\newblock {\em Illinois J. Math.}, 11:417--427, 1967.

\bibitem[Hop41]{MR4784}
Heinz Hopf.
\newblock \"{U}ber die {T}opologie der {G}ruppen-{M}annigfaltigkeiten und ihre
  {V}erallgemeinerungen.
\newblock {\em Ann. of Math. (2)}, 42:22--52, 1941.

\bibitem[IK07]{MR2272137}
Norio Iwase and Akira Kono.
\newblock Lusternik-{S}chnirelmann category of {${\rm Spin}(9)$}.
\newblock {\em Trans. Amer. Math. Soc.}, 359(4):1517--1526, 2007.

\bibitem[IKM16]{MR3520436}
Norio Iwase, Kai Kikuchi, and Toshiyuki Miyauchi.
\newblock On {L}usternik-{S}chnirelmann category of {${\bf SO}(10)$}.
\newblock {\em Fund. Math.}, 234(3):201--227, 2016.

\bibitem[IM89]{MR1000378}
Norio Iwase and Mamoru Mimura.
\newblock Higher homotopy associativity.
\newblock In {\em Algebraic topology ({A}rcata, {CA}, 1986)}, volume 1370 of
  {\em Lecture Notes in Math.}, pages 193--220. Springer, Berlin, 1989.

\bibitem[IM04]{MR2039767}
Norio Iwase and Mamoru Mimura.
\newblock L-{S} categories of simply-connected compact simple {L}ie groups of
  low rank.
\newblock In {\em Categorical decomposition techniques in algebraic topology
  ({I}sle of {S}kye, 2001)}, volume 215 of {\em Progr. Math.}, pages 199--212.
  Birkh\"{a}user, Basel, 2004.

\bibitem[IS10]{MR2556074}
Norio Iwase and Michihiro Sakai.
\newblock Topological complexity is a fibrewise {L}-{S} category.
\newblock {\em Topology Appl.}, 157(1):10--21, 2010.

\bibitem[IS12]{MR2923451}
Norio Iwase and Michihiro Sakai.
\newblock Erratum to ``{T}opological complexity is a fibrewise {L}-{S}
  category'' [{T}opology {A}ppl. 157 (1) (2010) 10--21] [mr2556074].
\newblock {\em Topology Appl.}, 159(10-11):2810--2813, 2012.

\bibitem[IST19]{MR3975098}
Norio Iwase, Michihiro Sakai, and Mitsunobu Tsutaya.
\newblock A short proof for {${\rm tc}(K)=4$}.
\newblock {\em Topology Appl.}, 264:167--174, 2019.

\bibitem[Iwa98]{MR1642747}
Norio Iwase.
\newblock Ganea's conjecture on {L}usternik-{S}chnirelmann category.
\newblock {\em Bull. London Math. Soc.}, 30(6):623--634, 1998.

\bibitem[Iwa01]{MR1808218}
Norio Iwase.
\newblock Co-{$H$}-spaces and the {G}anea conjecture.
\newblock {\em Topology}, 40(2):223--234, 2001.

\bibitem[Iwa02]{MR1905835}
Norio Iwase.
\newblock {$A_\infty$}-method in {L}usternik-{S}chnirelmann category.
\newblock {\em Topology}, 41(4):695--723, 2002.

\bibitem[Iwa07]{MR2331629}
Norio Iwase.
\newblock The {G}anea conjecture and recent developments on
  {L}usternik-{S}chnirelmann category [translation of {S}\={u}gaku {\bf 56}
  (2004), no. 3, 281--296; mr2086116].
\newblock volume~20, pages 43--63. 2007.
\newblock Sugaku Expositions.

\bibitem[Iwa09]{MR2503529}
Norio Iwase.
\newblock Categorical length, relative {L}-{S} category and higher {H}opf
  invariants.
\newblock In {\em Algebraic topology---old and new}, volume~85 of {\em Banach
  Center Publ.}, pages 205--224. Polish Acad. Sci. Inst. Math., Warsaw, 2009.

\bibitem[Iwa12]{AX12115741}
Norio Iwase.
\newblock {A}ssociahedra, {M}ultiplihedra and units in ${A}_{\infty}$ form, 11
  2012.
\newblock arXiv: https://arxiv.org/pdf/1211.5741.pdf.

\bibitem[Jam78]{MR516214}
I.~M. James.
\newblock On category, in the sense of {L}usternik-{S}chnirelmann.
\newblock {\em Topology}, 17(4):331--348, 1978.

\bibitem[JM91]{MR1130605}
I.~M. James and J.~R. Morris.
\newblock Fibrewise category.
\newblock {\em Proc. Roy. Soc. Edinburgh Sect. A}, 119(1-2):177--190, 1991.

\bibitem[LS34]{LS1934}
L.~Lusternik and L.~Schnirelmann.
\newblock {\em M{\'e}thodes topologiques dans les probl{\`e}mes variationnels}.
\newblock Number~1 in Actualit{\'e}s Scientifiques Industrielles: Expos{\'e}s
  sur l'analyse math. et ses applications. Hermann, 1934.

\bibitem[Mat76]{MR402694}
Michael Mather.
\newblock Pull-backs in homotopy theory.
\newblock {\em Canadian J. Math.}, 28(2):225--263, 1976.

\bibitem[Mil56a]{MR77122}
John Milnor.
\newblock Construction of universal bundles. {I}.
\newblock {\em Ann. of Math. (2)}, 63:272--284, 1956.

\bibitem[Mil56b]{MR77932}
John Milnor.
\newblock Construction of universal bundles. {II}.
\newblock {\em Ann. of Math. (2)}, 63:430--436, 1956.

\bibitem[Mil59]{MR100267}
John Milnor.
\newblock On spaces having the homotopy type of a {${\rm CW}$}-complex.
\newblock {\em Trans. Amer. Math. Soc.}, 90:272--280, 1959.

\bibitem[Mil67]{MR208595}
R.~James Milgram.
\newblock The bar construction and abelian {$H$}-spaces.
\newblock {\em Illinois J. Math.}, 11:242--250, 1967.

\bibitem[RO99]{MR1686579}
Yuli~B. Rudyak and John Oprea.
\newblock On the {L}usternik-{S}chnirelmann category of symplectic manifolds
  and the {A}rnold conjecture.
\newblock {\em Math. Z.}, 230(4):673--678, 1999.

\bibitem[Rud98]{MR1679849}
Yuli~B. Rudyak.
\newblock Category weight: new ideas concerning {L}usternik-{S}chnirelmann
  category.
\newblock In {\em Homotopy and geometry ({W}arsaw, 1997)}, volume~45 of {\em
  Banach Center Publ.}, pages 47--61. Polish Acad. Sci. Inst. Math., Warsaw,
  1998.

\bibitem[Rud99]{MR1686583}
Yuli~B. Rudyak.
\newblock On analytical applications of stable homotopy (the {A}rnold
  conjecture, critical points).
\newblock {\em Math. Z.}, 230(4):659--672, 1999.

\bibitem[Sak10]{MR2665235}
Michihiro Sakai.
\newblock {$A_\infty$}-spaces and {L}-{S} category in the category of fibrewise
  spaces.
\newblock {\em Topology Appl.}, 157(13):2131--2135, 2010.

\bibitem[Sin75]{MR391075}
Wilhelm Singhof.
\newblock On the {L}usternik-{S}chnirelmann category of {L}ie groups.
\newblock {\em Math. Z.}, 145(2):111--116, 1975.

\bibitem[Sin76]{MR425952}
Wilhelm Singhof.
\newblock On the {L}usternik-{S}chnirelmann category of {L}ie groups. {II}.
\newblock {\em Math. Z.}, 151(2):143--148, 1976.

\bibitem[Sta63]{MR0158400}
James~Dillon Stasheff.
\newblock Homotopy associativity of {$H$}-spaces. {I}, {II}.
\newblock {\em Trans. Amer. Math. Soc. 108 (1963), 275-292; ibid.},
  108:293--312, 1963.

\bibitem[Sta70]{MR0270372}
James Stasheff.
\newblock {\em {$H$}-spaces from a homotopy point of view}.
\newblock Lecture Notes in Mathematics, Vol. 161. Springer-Verlag, Berlin-New
  York, 1970.

\bibitem[Str97]{MR2696593}
Jeffrey~Andrew Strom.
\newblock {\em Category weight and essential category weight}.
\newblock ProQuest LLC, Ann Arbor, MI, 1997.
\newblock Thesis (Ph.D.)--The University of Wisconsin - Madison.

\bibitem[Sug57a]{MR97066}
Masahiro Sugawara.
\newblock A condition that a space is group-like.
\newblock {\em Math. J. Okayama Univ.}, 7:123--149, 1957.

\bibitem[Sug57b]{MR86303}
Masahiro Sugawara.
\newblock On a condition that a space is an {$H$}-space.
\newblock {\em Math. J. Okayama Univ.}, 6:109--129, 1957.

\bibitem[Tod62]{MR0143217}
Hirosi Toda.
\newblock {\em Composition methods in homotopy groups of spheres}.
\newblock Annals of Mathematics Studies, No. 49. Princeton University Press,
  Princeton, N.J., 1962.

\bibitem[\v{S}61]{MR0154284}
A.~S. \v{S}varc.
\newblock The genus of a fibered space.
\newblock {\em Trudy Moskov. Mat. Ob\v{s}\v{c}.}, 10:217--272, 1961.

\bibitem[\v{S}62]{MR0151982}
A.~S. \v{S}varc.
\newblock The genus of a fibre space.
\newblock {\em Trudy Moskov Mat. Ob\v{s}\v{c}.}, 11:99--126, 1962.

\bibitem[\v{S}66]{AMST1966Svarc}
A.~S. \v{S}varc.
\newblock The genus of a fiber space.
\newblock {\em American Mathematical Society Translations: Series 2},
  55:49--140, 1966.

\end{thebibliography}

\end{document}